\documentclass[leqno]{siamltex}

\title{Risk-Sensitive Markov Control Processes}

\author{Yun Shen\thanks{Fakult\"at Elektrotechnik und Informatik, Technische Universit\"at Berlin,  Marchstr.~23, 10587, Berlin, Germany ({\tt yun@ni.tu-berlin.de}).} \and Wilhelm Stannat\thanks{Institut f\"ur Mathematik, and Bernstein Center for Computational Neuroscience, Technische Universit\"at Berlin, Stra{\ss}e des 17.~Juni 136, 10623, Berlin, Germany ({\tt stannat@math.tu-berlin.de}). The work of this author was supported by the BMBF, FKZ01GQ1001B.}
 \and Klaus Obermayer\thanks{Fakult\"at Elektrotechnik und Informatik, and Bernstein Center for Computational Neuroscience, Technische Universit\"at Berlin, Marchstr.~23, 10587, Berlin, Germany ({\tt oby@ni.tu-berlin.de}). The work of this author was supported by the BMBF (Bersteinfokus Lernen TP1), 01GQ0911.}}

\usepackage{amsmath}
\usepackage{amsfonts}
\usepackage{mathrsfs}
\usepackage{enumerate}
\usepackage{url}

\newtheorem{assumption}{Assumption}[section]

\newcommand{\maps}{\rightarrow}

\newcommand{\diff}{\textrm{d}}
\newcommand{\bs}{\boldsymbol}

\newcommand{\ftk}{risk map}
\newcommand{\fpm}{risk measure}
\newcommand{\ftkc}{Risk Map}
\newcommand{\fpmc}{Risk Measure}

\begin{document}
 \maketitle

\begin{abstract}
We introduce a general framework for measuring risk in the context of Markov control processes with risk maps on general Borel spaces that generalize known concepts of risk measures in mathematical finance, operations research and behavioral economics. Within the framework, applying weighted norm spaces to incorporate also unbounded costs, we study two types of infinite-horizon risk-sensitive criteria, discounted total risk and average risk, and solve the associated optimization problems by dynamic programming. For the discounted case, we propose a new discount scheme, which is different from the conventional form but consistent with the existing literature, while for the average risk criterion, we state Lyapunov-like stability conditions that generalize known conditions for Markov chains to ensure the existence of solutions to the optimality equation.
\end{abstract}



\pagestyle{myheadings}
\thispagestyle{plain}

\section{Introduction}

\emph{Markov control processes} (MCPs, see
e.g.~\cite{hernandez1996discrete,hernandez1999further} and
\cite{puterman1994markov} under the name \emph{Markov decision processes})
are widely applied to model sequential decision making problems of
agents. The induced optimal control problem is to find the best policy that
minimizes the expected total costs. The core of the MCP-framework consists of two \emph{objective} descriptions of some mechanism of the environment: \emph{transition probabilities} of switching states when
performing actions, and immediate \emph{outcomes} (rewards or costs)
obtained at states by executing actions.  Facing the same environment, however, different agents might have
different policies. Therefore, in many applications, it is important to also
incorporate the \emph{subjective} perceptions of an agent into the
MCP-framework. The subjective outcomes are usually modeled by \emph{utility
  functions} (see e.g.~\cite{gollier2004economics}), which can be easily incorporated by simply
replacing the immediate outcome with its utility, whereas the subjective
transition probabilities require a more sophisticated mathematical
framework. They are commonly incorporated in the \emph{risk}, which is caused by an
uncertain environment.

\emph{Coherent/convex risk measures} (CRMs) \cite{artzner1999coherent,follmer2002convex} have been widely employed to model subjective probabilities in mathematical finance since the last decade. Several works (see e.g.~\cite{roorda2005coherent,foellmer2006convex,cheridito2006dynamic,ruszczynski2006conditional,cheridito2011composition} and references therein) extend CRMs to temporal structures in various setups, where they consider mainly finite-horizon problems. On the contrary, in the literature of MCPs, while the infinite-horizon risk-sensitive optimal control problems are studied, they apply merely the \emph{entropic map} \cite{chung1987discounted, avila1998controlled, borkar2002risk, di2008infinite,coraluppi2000mixed,fleming1997risk,hernandez1996risk,
  marcus1997risk, cavazos2010optimality}, which is convex and in fact a special type of CRM. All risk measures mentioned in the above literature are coherent/convex based on the assumption that the agent is supposed to be economically rational and therefore \emph{risk-averse}. This limits applications in the fields of decision-making under risk and behavioral economics, where more general risk measures (see e.g.~\cite{savage1972foundations, chateauneuf2008cardinal, tversky1992advances} and references therein) are applied, since  human beings are not always risk-averse. However, the models in these fields can only be applied to one-step decision making problems.

To overcome the limitations mentioned above, (i) we extend the definition of CRMs to include the risk measures considered also in behavioral economics; (ii) we apply a \emph{constructive approach} (see \eqref{eq:rdcmp1} below) which maintains the Markov property that is necessary for the existence of stationary optimal policies for two infinite-horizon objectives, albeit less general than the risk maps used in \cite{ruszczynski2006conditional,ruszczynski2010risk}.

More specifically, three types of objectives are usually considered in the literature of MCPs: finite-stage, discounted and average cost, depicted as
\begin{equation}
 S_T := \sum_{t=0}^T c(X_t,A_t), \quad S_\alpha := \sum_{t=0}^\infty \alpha^t c(X_t,A_t), \quad \textrm{and } \quad S_A := \limsup_{T \rightarrow \infty} \dfrac{1}{T} S_T,  \label{eq:totalreward1}
\end{equation}
where $X_t$ and $A_t$ are state and action at time $t$ respectively, $c$ denotes the cost function and $\alpha \in [0,1)$ is the discount factor. Given the initial state $X_0 = x$, the optimization problem is then to minimize the expected objective
\begin{equation}
 \inf_{\bs \pi} \mathbb E^{\bs \pi} \left[ \mathcal S | X_0 = x\right] \label{eq:obj1}
\end{equation}
by selecting a Markov policy $\bs \pi = [\pi_0,\pi_1,\ldots]$, where $\mathcal S$ is $S_T$, $S_\alpha$ or $S_A$. We are mainly interested in the risk-sensitive extensions of two infinite-horizon objectives, i.e., the discounted and average one. We notice that the discounted objective (while other two objectives can be dealt with analogously) can be decomposed as follows,
\begin{equation}
 \mathbb E^{\bs \pi}_{X_0} \left[ S_\alpha  \right] = c^{\pi_0}(X_0) +  \alpha \mathbb E^{\pi_0}_{X_0} \left[
   c^{\pi_1}(X_1) + \alpha \mathbb E^{\pi_1}_{X_1} \left[ c^{\pi_2}(X_2) + \ldots \right]  \right] \label{eq:decompose}
\end{equation}
where $\mathbb E^{\pi}_{X}[\cdot]$ denotes conditional expectation at state $X$ under policy $\pi$. We replace the risk-neutral conditional expectation with \emph{risk maps} $\mathcal R^{\pi}_{X}[\cdot]$ similar to the Markov risk measures defined in \cite{ruszczynski2006conditional,ruszczynski2010risk} and obtain the \emph{risk-sensitive} objective
\begin{equation}
J^{\bs \pi}_\alpha = c^{\pi_0}(X_0) +  \alpha \mathcal R^{\pi_0}_{X_0} \left[
   c^{\pi_1}(X_1) +  \alpha  \mathcal R^{\pi_1}_{X_1} \left[ c^{\pi_2}(X_2) + \ldots  \right]
 \right]. \label{eq:rdcmp1}
\end{equation}

With the generalized risk measures and constructed risk maps, we provide a unified treatment in the context of MCPs to infinite-horizon risk-sensitive optimal control problems considered in various fields,  e.g.~optimal control, operations research, finance and behavioral economics. Using weighted norm spaces, we can for the first time incorporate unbounded costs in risk-sensitive MCPs also. We prove that two types of objectives, the discounted total risk and the average risk, can be optimized with \emph{dynamic programming} algorithms under proper assumptions. For the case of discounted risk, we apply a new discount scheme which is different from the conventional form but consistent with the one applied in \cite{ruszczynski2010risk} where coherent risk measures are considered. For the average case, we state sufficient conditions, which generalize Lyapunov-type conditions from the literature of Markov chains (see e.g.~\cite{meyn1993markov}), to ensure the existence of solutions to the associated optimality equation.

The paper is organized as follows. In Section \ref{sec:preli} we introduce our concept of risk measures in the context of MCPs on Borel spaces, generalizing CRMs considered in the mathematical finance, in order to also include the family of risk measures considered in behavioral economics. In Section \ref{sec:wc}, we extend our definition of risk measures to a Markovian temporal structure and call them accordingly \emph{risk maps}, whose properties are also investigated. In Section \ref{sec:armtmcp}, we consider risk maps in the MCP-framework by adding control parameters. We demonstrate in Subsection \ref{sec:pre} how to explicitly construct complex risk maps by combining simpler ones, e.g., risk-averse and risk-seeking maps, followed by a discussion of examples covered by our framework in Subsection \ref{sec:ex}. The induced infinite-horizon risk-sensitive objectives, including the discounted and average risk, are optimized under assumptions in Subsection \ref{sec:disc} and \ref{sec:ar}. Finally, in Section \ref{sec:rm}, we present one example with one risk map and prove that the proposed sufficient conditions for average risk are satisfied.

\section{Preliminaries}
\label{sec:preli}
In this section, we introduce the framework of Markov control processes (MCPs) and \fpm{s}, where in the first part, we mostly follow the notations of Hern{\'a}ndez-Lerma \& Lasserre (1999) \cite{hernandez1999further}. We clarify some concepts before stating setups. A \emph{Borel space} is a Borel subset of a complete separable metric space. If $\mathbf X$ is a Borel space, its Borel $\sigma$-algebra is denoted by $\mathcal B(\mathbf X)$. Let $\mathbf X$ and $\mathbf Y$ be two Borel spaces. A \emph{stochastic kernel on $\mathbf X$ given $\mathbf Y$} is a function $\psi(B|y), B \in \mathcal B(\mathbf X), y \in \mathbf Y$ such that i) $\psi(\cdot|y)$ is a probability measure on $\mathcal B(\mathbf X)$ for every fixed $y \in \mathbf Y$, and ii) $\psi(B|\cdot)$ is a measurable function on $\mathbf Y$ for every fixed $B \in \mathcal B(\mathbf X)$. 

\subsection{Markov Control Processes}
\label{sec:mcp}
A Markov control process, $(\mathbf X, \mathbf A,\{\mathbf A(x)|x \in \mathbf X \},Q,c)$, consists of  the following components: \emph{state space} $\mathbf X$ and \emph{action space} $\mathbf A$, which are Borel spaces; the feasible action set $\mathbf A(x)$, which is a nonempty Borel space of $\mathbf A$, for a given state $x \in \mathbf X$; the \emph{transition model} $Q(B|x,a), B \in \mathcal B(\mathbf X), (x,a) \in \mathbf K$: a \emph{stochastic kernel} on $\mathbf X$ given $\mathbf K$, where $\mathbf K$ denotes the set of feasible state-action pairs $\mathbf K:=\{ (x,a)|x \in \mathbf X, a \in \mathbf A(x)\}$, which is a Borel subset of $\mathbf X \times \mathbf A$; and the \emph{cost function} $c$: $\mathbf K \maps \mathbb R$, $\mathcal B(\mathbf K)$-measurable. Random variables are denoted by capital letters, e.g.~$X_t$ and $A_t$, whereas realizations of the random variables are denoted by normal letters, e.g.~$x_t$ and $a_t$.

We consider in this paper \emph{Markov policies}, $\bs \pi = [ \pi_0,\pi_1,\pi_2,\ldots]$, where each \emph{single-step policy} $\pi_t(\cdot |x_t)$, which denotes the probability of choosing action $a_t$ at $x_t$, $(x_t,a_t) \in \mathbf K$, is Markov (independent of the states and actions before $t$) and, therefore, a stochastic kernel on $\mathbf A$ given $\mathbf X$. We use the bold typeface to represent a sequence of policies while using normal typeface for a single-step policy. Let $\Delta$ denote the set of all stochastic kernels on $\mathbf A$ given $\mathbf X$, $\mu$, such that $\mu(\mathbf A(x)|x) = 1$ and $\Pi_M$ denotes the set of all Markov policies. Thus $\Pi_M = \Delta^\infty$. A policy $f \in \Delta$ is \emph{deterministic} if for each $x\in \mathbf X$, there exists some $a\in \mathbf A(x)$ such that $f(\{a\}|x) = 1$. Let $\Delta_D \subset \Delta$ denote the set of all deterministic single-step policies. A policy $\bs \pi $ is said to be \emph{stationary}, if $\bs \pi = \pi^\infty$ for some $\pi \in \Delta$. For each $x \in \mathbf X$ and single-step policy $\pi \in \Delta$, define
\begin{align}
 c^\pi(x) := \int_{\mathbf A(x)} c(x,a) \pi(\diff a|x), P^\pi(B|x) := \int_{\mathbf A(x)} Q(B|x,a) \pi(\diff a|x), B \in \mathcal B(\mathbf X). \label{eq:cpi}
\end{align}

There are usually three types of objectives used in the literature of MCPs: finite-stage, discounted and average cost, depicted as
\begin{equation}
 S_T := \sum_{t=0}^T c(X_t,A_t), \quad S_\alpha := \sum_{t=0}^\infty \alpha^t c(X_t,A_t), \quad \textrm{and } \quad S_A := \limsup_{T \rightarrow \infty} \dfrac{1}{T} S_T  \label{eq:totalreward}
\end{equation}
where $\alpha \in [0,1)$ denotes the discount factor. Suppose we start from one given state $X_0 = x$. The optimization problem is then to minimize the expected objective
\begin{equation}
 \inf_{\bs \pi \in \Pi_M} \mathbb E^{\bs \pi} \left[ \mathcal S | X_0 = x\right] \label{eq:obj}
\end{equation}
by selecting a policy $\bs \pi$, where $\mathcal S$ is $S_T$, $S_\alpha$ or
$S_A$. We notice that the finite-stage objective function (while other two objectives can be dealt with analogously) can be decomposed as follows,
\begin{equation*}
 \mathbb E^{\bs \pi}_{X_0} \left[ S_T  \right] = c^{\pi_0}(X_0) +  \mathbb E^{\pi_0}_{X_0} \left[
   c^{\pi_1}(X_1) + \mathbb E^{\pi_1}_{X_1} \left[ c^{\pi_2}(X_2) + \ldots + \mathbb
     E^{\pi_{T-1}}_{X_{T-1}} \left[ c^{\pi_T}(X_T) \right] \ldots  \right]
 \right] 
\end{equation*}
where $ \mathbb E^{\pi_t}_{X_t}\left[ v(X_{t+1}) \right] := \int v(X_{t+1}) P^{\pi_t}(dX_{t+1} |X_t)$ denotes the \emph{conditional expectation} of the function $v$ of the successive state $X_{t+1}$ given current state $X_t$. Obviously, the conditional expectation plays the key role in the calculation of all three objectives. 

In order to incorporate risk, we directly replace the expectation $\mathbb E_{X_t}^{\pi_t}$ with a \emph{\ftk{}} $\mathcal R_{X_t}^{\pi_t}$ that is similar to the risk mapping defined in \cite{ruszczynski2006conditional} and will be formally introduced in Section \ref{sec:armtmcp}. With the replacement, we obtain the \emph{risk-sensitive} objective
\begin{equation}
J^{\bs \pi}_T = c^{\pi_0}(X_0) +  \mathcal R^{\pi_0}_{X_0} \left[
   c^{\pi_1}(X_1) + \mathcal R^{\pi_1}_{X_1} \left[ c^{\pi_2}(X_2) + \ldots + \mathcal R^{\pi_{T-1}}_{X_{T-1}} \left[ c^{\pi_T}(X_T) \right] \ldots  \right]
 \right]. \label{eq:rdcmp}
\end{equation}
The other two objectives will be defined analogously and discussed in Section \ref{sec:optrisk}.

{\em Remark.} Comparing with \cite{ruszczynski2010risk} where dynamic risk measures depending on the whole history are allowed, we apply here a less general and Markovian type of risk measures which depends only on the current state. It is due to the motivation that (i) the underlying structure is Markovian, and (ii) for optimizing infinite-horizon objectives, it becomes computationally infeasible if risk measures are dependent on the whole history.

\subsection{\fpmc{s}}
\label{sec:fpm}
In order to include risk measures considered also in behavioral economics, we introduce a generalized version of risk measures which are originally defined in \cite{artzner1999coherent,follmer2002convex}. Consider a probability space $(\mathbf X, \mathcal B(\mathbf X))$ under some probability measure $\varphi$. Let $\mathscr L$ be a vector space of $\mathcal B(\mathbf X)$-measurable real-valued functions and $\mathscr B$ be the space of all bounded $\mathcal B(\mathbf X)$-measurable real-valued functions. We assume $\mathscr L \supset \mathscr B$.  The partial ordering ``$\leq$'' in $\mathscr L$ is defined as: $v \leq u$, if $v(x) \leq u(x)$ for all $x \in \mathbf X$. For convenience, we write $u \in \mathbb R$ if $u$ is a constant real-valued function. Thus $u \in \mathscr L$, which will be specified in Section \ref{sec:wn}.
\begin{definition} A mapping $\nu: \mathscr L \maps \mathbb R \cup \{ \infty\}$ is said to be a \textbf{\fpm} if 
 \begin{enumerate}[(I)]
  \item (Monotonicity) $\nu(v) \leq \nu(u)$, whenever $v\leq u \in \mathscr L$;
  \item (Translation invariance) $\nu(v + u) = \nu(v) + u$, for any $u \in \mathbb R$;
  \item (Centralization) $\nu(0) = 0$.
 \end{enumerate}
Moreover, $\nu$ is called real-valued if $\nu(v) \in \mathbb R$ for all $v \in \mathscr L$. 
\label{def:rm}
\end{definition}

{\em Remark.} Within the economic context, $v$ and $u$ are usually considered as random variables, which are used to model the \emph{uncertain} costs in the future. Monotonicity reflects the intuition that if the cost of one case is higher than the cost of another case, the risk of the case must be higher than that of the other one. Under the axiom of translation invariance, the sure cost $u$ (equal everywhere in the state space $\mathbf X$) in the future, which can be viewed as a constant function, is considered as a sure cost at current time point. The axiom of centralization sets the reference point to be 0, i.e., there is no risk if there is no cost.

\begin{definition}\label{def:cat} A \fpm{} $\nu$ is called 
\begin{itemize}
 \item \textbf{convex}, if for all $\alpha \in [0,1]$, $v, u \in \mathscr L$, $\nu(\alpha v+ (1-\alpha)u) \leq  \alpha \nu( v) + (1-\alpha)\nu(u)$;
 \item \textbf{concave}, if for all $\alpha \in [0,1]$, $v, u \in \mathscr L$, $\nu(\alpha v+ (1-\alpha)u) \geq  \alpha \nu( v) + (1-\alpha)\nu(u)$;
 \item \textbf{homogeneous}, if for all $\lambda \in \mathbb R_+$ and $v \in \mathscr L$, $\nu(\lambda v ) =   \lambda \nu(v)$;
 \item \textbf{coherent}, if $\nu$ is convex and homogeneous.
\end{itemize}
\end{definition}

{\em Remark.} Comparing with the risk measures defined in finance \cite{artzner1999coherent,follmer2002convex}, neither convexity nor coherency is required. We will see in Section \ref{sec:pre} that convex \fpm{s} correspond to the case that
the agent is risk-averse. However, in some problems, especially in modeling real human behaviors, mixed risk-preference (risk-averse at some states while risk-seeking at other states) is also a possible strategy. For instance, at gambling, some people are risk-averse when losing money but risk-seeking when winning money. Therefore, we require neither coherence nor convexity. 

\section{\ftkc{s} without Control}
\label{sec:wc}
Before applying \ftk{s} in the MCP-framework, we first define \ftk{s} without control on Markov chains and investigate their properties. We refer to \cite{ruszczynski2006conditional} for conditional risk maps under more general settings. Consider a time-homogenous Markov chain with state space $\mathbf X$ and transition kernel $P$. 
\begin{definition}\label{def:rmaps}
A mapping $\mathcal R(x,v): \mathbf X \times \mathscr L \maps \mathbb R \cup \{ \infty \}$ is said to be a \textbf{\ftk} on the Markov chain $P$, if (i) for each $x \in \mathbf X$, $\mathcal R_x(\cdot) := \mathcal R(x,\cdot) $ is a \fpm{}; and (ii) $\mathcal R(\cdot, v) \in \mathscr L$ for each $v \in \mathscr L$. $\mathcal R$ is called real-valued if for each $x \in \mathbf X$ and $v \in \mathscr L$, $\mathcal R(x,v) \in \mathbb R$. $\mathcal R$ is called convex (resp.\ concave, homogenous, coherent) if for all $x \in \mathbf X$, $\mathcal R_x$ is convex (resp.\ concave, homogenous, coherent). 
\end{definition}

{\em Remark.} With slight abuse of terminology, $P$ can be viewed as a linear operator such that $P_x (v) := \int v(y) P(dy|x)$, $v \in \mathscr L$. Then the operator $\mathcal R$ defined above is a generalization of $P$, maintaining the key properties, monotonicity, translation invariance and centralization. In the following, we view $\mathcal R(v)$ as a function on $\mathscr L$, for each $v \in \mathscr L$. Finally, for two risk maps, $\mathcal S$, $\mathcal R$, we write $\mathcal S \leq \mathcal R$, if $\mathcal S(v) \leq \mathcal R(v)$ for all $v \in \mathscr L$.
\subsection{Sub- and Uppermodules}
\label{sec:mod}
To investigate stability properties of risk maps, we generalize the modules introduced by Delbaen (2000) \cite{Delbaen_2000}. 
\begin{definition}
Let $\mathcal F(x,v): \mathbf X \times \mathscr L \maps \mathbb R $ be a map on $\mathbf X$ and $\mathscr L$ and write $\mathcal F_x(v) := \mathcal F(x,v)$. Then $\mathcal F^\sharp: \mathbf X \times \mathscr L \maps \mathbb R \cup\{\infty\}$ defined by $\mathcal F^\sharp_x(v) := \sup_{u \in \mathscr L} \left\{ \mathcal F_x(v + u) - \mathcal F_x(u)\right\}$, $\forall x \in \mathbf X,$ is called the \textbf{submodule} of $\mathcal F$.   
Let $\mathcal G(x,v): \mathbf X \times \mathscr L \maps \mathbb R \cup \{ \infty \}$ be a map on $\mathbf X$ and $\mathscr L$ and write $\mathcal G_x(v) := \mathcal G(x,v)$. Then $\overline{\mathcal G}: \mathbf X \times \mathscr L \maps \mathbb R \cup\{\infty\}$ defined by $ \overline{\mathcal G}_x(v) := \sup_{\lambda \neq 0} \frac{\mathcal G_x(\lambda v)}{\lambda}$, $\forall x \in \mathbf X$, is called the \textbf{uppermodule} of $\mathcal G$. 
\end{definition}

We summarize some properties of sub- and uppermodules of risk maps.
\begin{proposition}\label{prop:sm}
  Let $\mathcal R$ be a real-valued \ftk{}. Then for each $x \in \mathbf X$ (i) $\mathcal R_x(v) \leq \mathcal R^\sharp_x(v)$, $\forall v \in \mathscr L$; (ii) $\mathcal R^\sharp_x$ is a risk measure; (iii) $\mathcal R^\sharp_x$ is sublinear, i.e., $\mathcal R^\sharp_x(v + v') \leq \mathcal R^\sharp_x(v) + \mathcal R^\sharp_x(v')$.
\end{proposition}
\begin{proof}
(i) and (ii) are immediate consequences of the construction of the submodule $\mathcal R^\sharp$ associated with $\mathcal R$. It remains to prove (iii). First, for all $v \in \mathscr L$ and $x \in \mathbf X$, $\mathcal R^\sharp_x(v) = \sup_{u \in \mathscr L} \left\{ \mathcal R_x(v + u) - \mathcal R_x(u)\right\} = \sup_{u' \in \mathscr L} \left\{ \mathcal R_x(u') - \mathcal R_x(u'-v)\right\}$. Thus, setting $u' = v' + u$, we have for all $u \in \mathscr L $ and $x \in \mathbf X$,
\begin{align*}
  &\mathcal R_x(v + v' + u) - \mathcal R_x(u) \leq \sup_{u' \in \mathscr L} \left\{ \mathcal R_x(v + u') - \mathcal R_x(u') + \mathcal R_x(u') - \mathcal R_x(u'-v') \right\} \\
 \leq & \sup_{u' \in \mathscr L} \left\{ \mathcal R_x(v + u') - \mathcal R_x(u') \right\} + \sup_{u' \in \mathscr L} \left\{ \mathcal R_x(u') - \mathcal R_x(u'-v') \right\} =  \mathcal R^\sharp_x(v) + \mathcal R^\sharp_x(v').
\end{align*}
Hence $\mathcal R^\sharp(v + v') \leq \mathcal R^\sharp(v) + \mathcal R^\sharp(v')$.
\quad \end{proof}

\begin{proposition}\label{prop:um}
 Let $\mathcal R(x,v)$ be a map on $\mathbf X \times \mathscr L$ such that for each $x \in \mathbf{X}$, $\mathcal R(x,\cdot)$ is a risk measure. Then for all $x\in \mathbf{X}$, (i) $\mathcal R_x(v) \leq \overline{\mathcal R}_x(v)$, $\forall v \in \mathscr L$; (ii) $\overline{\mathcal R}_x$ is a risk measure; (iii) $\overline{\mathcal R}_x$ is homogeneous.
\end{proposition}

{\em Proof.}
 Similar to their counterparts for submodule, (i) and (ii) are immediate consequences of the construction of the uppermodule $\overline{\mathcal R}$ associated with $\mathcal R$. For the proof of (iii), note that for any $\beta > 0$,
$$
 \overline{\mathcal R}_x(\beta v) = \sup_{\lambda \neq 0} \dfrac{\mathcal R_x(\lambda \beta v)}{\lambda} = \sup_{\lambda \beta \neq 0} \beta \dfrac{\mathcal R_x(\lambda \beta v)}{\beta \lambda} = \beta \sup_{\lambda \neq 0} \dfrac{\mathcal R_x(\lambda v)}{\lambda} = \beta \overline{\mathcal R}_x(v).  \eqno \endproof
$$

\begin{proposition}\label{prop:co}
 Let $\mathcal R$ be a real-valued \ftk. Then $\forall x \in \mathbf X$, $\overline{\mathcal R^\sharp}_x(v) = \sup_{\lambda > 0} \frac{\mathcal R_x^\sharp(\lambda v)}{\lambda}$ and $\overline{\mathcal R^\sharp}_x$ is coherent. Furthermore, if $\mathcal R$ is coherent, then $\overline{\mathcal R^\sharp} = \mathcal R$.
\end{proposition}
\begin{proof}
 (i) First we show that $- \lambda^{-1}\mathcal R^\sharp_x( -\lambda v ) \leq \lambda^{-1} \mathcal R^\sharp_x( \lambda v )$ for each $\lambda \in \mathbb R_+$ and $v \in \mathscr L$. Indeed, without loss of generality, assume that $R^\sharp_x( - \lambda v )< \infty$. Then by sublinearity, $0 = \mathcal R^\sharp_x(0) = \mathcal R^\sharp_x(\lambda v - \lambda v)  \leq \mathcal R^\sharp_x(\lambda v )+\mathcal R^\sharp_x( -\lambda v )$ implies that $- \lambda^{-1}\mathcal R^\sharp_x( -\lambda v ) \leq \lambda^{-1} \mathcal R^\sharp_x( \lambda v )$. Thus $\overline{\mathcal R^\sharp}_x(v) = \sup_{\lambda > 0} \frac{\mathcal R_x^\sharp(\lambda v)}{\lambda}, \forall x \in \mathbf X$. $\overline{\mathcal R^\sharp}_x$ is coherent, since $\overline{\mathcal R^\sharp}_x$ is both sublinear and homogeneous. 

(ii) Suppose $\mathcal R$ is coherent. Then $\mathcal R(v + u ) - \mathcal R(u ) \leq \mathcal R(v)$ for all $v, u \in \mathscr L$. Hence, $\mathcal R^\sharp \leq \mathcal R$. Together with Proposition \ref{prop:sm}(i) $\mathcal R \leq \mathcal R^\sharp$, we have $\mathcal R^\sharp = \mathcal R$. Since $\mathcal R$ is homogenous, $\overline{\mathcal R^\sharp} = \mathcal R$. \quad \end{proof}

\begin{proposition}\label{prop:sharptop}
Let $\mathcal R$ be a real-valued risk map. Then $\lvert \mathcal R_x(v + u) - \mathcal
 R_x(u) \rvert \leq \overline{\mathcal R^\sharp}_x(\lvert v \rvert)$ for all $v, u \in \mathscr L$ and $x \in \mathbf X$.
\end{proposition}
\begin{proof}
By Proposition \ref{prop:sm} and \ref{prop:um}, for each $x\in \mathbf X$, $ \mathcal R_x(v + u) - \mathcal
 R_x(u) \leq \overline{\mathcal R^\sharp}_x(v) \leq \overline{\mathcal R^\sharp}_x(\lvert v \rvert)$. On the other hand, $\mathcal R_x(u) - \mathcal  R_x(u+v) \leq \overline{\mathcal R^\sharp}_x(-v) \leq \overline{\mathcal R^\sharp}_x(\lvert v \rvert)$.
\quad \end{proof}

{\em Remark.} For any real-valued risk map $\mathcal R$, its sub-upper-module $\overline{\mathcal R^\sharp}$ satisfying the three axioms of risk measures is in fact its coherent upper bound. In following sections, we shall apply this upper bound to control the growth speed of iterations of $\mathcal R$ (see Assumption \ref{asmp:suff}).

\subsection{Weighted Norm}
\label{sec:wn}
We now specify the functional space $\mathscr L$ used in this paper. Suppose $w: \mathbf X \maps [1,\infty)$ is a given measurable function. Consider the $w$-norm
\begin{align*}
 \lVert u \rVert_{w} := \sup_{x \in \mathbf X} \dfrac{\lvert u(x) \rvert }{w(x)}.
\end{align*}
Let $\mathscr B_w$ be the space of real-valued $w$-bounded $\mathcal B(\mathbf X)$-measurable functions. It is obvious that the bounded functional space $\mathscr B \subset \mathscr B_w$. Let $\mu$ be a signed measure on $\mathcal B(\mathbf X)$. Define $\lVert \mu \rVert_w := \sup_{\lVert u \rVert_w \leq 1} \lvert \int_{\mathbf X} u d \mu \rvert = \int_{\mathbf X} w \diff \lvert \mu \rvert \geq \lVert \mu \rVert_{TV}$. Denote by $\mathscr M$ (resp.~$\mathscr M_w$) the space of all (resp.~$w$-)bounded signed measures on $\mathcal B(\mathbf X)$. Thus $\mathscr M_w \subset \mathscr M$. More discussions about the $w$-norm space with applications in MCPs are referred to \cite{hernandez1999further}.

\begin{proposition}\label{prop:wbar}
 Let $\mathcal R$ be a real-valued risk map. Suppose there exists a $\overline w \in \mathbb R_+$ such that $\overline{\mathcal R^\sharp}_x(w) \leq \overline w \cdot w(x), \forall x \in \mathbf X$. Then $\lVert \mathcal R(v)- \mathcal R(u) \rVert_{w} \leq \overline w \lVert v - u \rVert_{w}$, for all $v,u \in \mathscr B_w$.
\end{proposition}
\begin{proof}
By definition, $\lvert v - u\rvert \leq \lVert v - u \rVert_{w} w$. Hence, Proposition \ref{prop:sharptop} yields
\begin{align*}
 \lvert \mathcal R(v) - \mathcal R (u) \rvert \leq \overline{\mathcal R^\sharp}(\lvert v - u \rvert) \leq  \overline{\mathcal R^\sharp}(\lVert v - u \rVert_w w) = \lVert v - u \rVert_w \overline{\mathcal R^\sharp}(w)
\end{align*}
due to the homogeneity of $\overline{\mathcal R^\sharp}$. Using the assumption $\overline{\mathcal R^\sharp}(w) \leq \overline w \cdot w$, we obtain that $\lVert \mathcal R(v)- \mathcal R(u) \rVert_{w} \leq \overline w \lVert v - u \rVert_{w}$.
\quad \end{proof}

\begin{corollary}
 Under the same assumption of Proposition \ref{prop:wbar}, $\mathcal R(v) \in \mathscr B_{w}$, for all $v \in \mathscr B_{w}$.
\end{corollary}

We consider the following seminorm that will play a key role in Section \ref{sec:pe} to investigate the stability of \ftk{s} and the existence of a solution to the Poisson equation.
\begin{align*}
 \lVert v \rVert_{s,w} := \sup_{x \neq y} \dfrac{\lvert v(x) - v(y) \rvert}{d_w(x,y)}, \textrm{where } d_w(x,y) := \left\{ \begin{array}{ll}
                      0 & x=y \\
w(x) + w(y) & x\neq y
                     \end{array}
 \right. .
\end{align*}
This seminorm is originally used by Hairer and Mattingly (2011) \cite{hairer2011yet} to study the ergodicity of Markov chains. In particular, when restricting to the bounded space $\mathscr B$, i.e.~setting $w \equiv 1$, the seminorm is called \emph{span-norm} in \cite{hernandez1989adaptive} and \emph{Hilbert seminorm} in \cite{gaubert2004perron}. In the following, we restate the Lemma 2.1 in \cite{hairer2011yet} and incorporate its proof for readers' convenience.
\begin{lemma}\label{lm:semi}
$\lVert v \rVert_{s,w} = \inf_{c \in \mathbb R} \lVert v + c \rVert_w, \forall v \in \mathscr B_w$.
\end{lemma}
\begin{proof}
 It is obvious that $\lVert v \rVert_{s,w} \leq \lVert v \rVert_w$ and therefore $\lVert v \rVert_{s,w} \leq \inf_{c \in \mathbb R} \lVert v +c \rVert_{w}$. It remains to prove the reverse inequality. Given any $\lVert v \rVert_{s,w} \leq 1$, set $c = \inf_{x} \{ w(x) - v(x)\}$. Note that for any $x$ and $y$, $v(x) \leq \lvert v(y) \rvert + \lvert v(x) - v(y) \rvert \leq \lvert v(y) \rvert + w(x) + w(y)$. Hence $w(x) - v(x) \geq - w(y) - \lvert v(y) \rvert$, which implies that $c$ is bounded below and hence $\lvert c \rvert < \infty$. Observe that $v(x) + c \leq v(x) + w(x) - v(x) \leq w(x)$ and $$v(x) + c = \inf_y \{ v(x) + w(y) - v(y)\} \geq \inf_y \{ w(y) - d_w(x,y) \lVert v \rVert_{s,w} \} \geq - w(x).$$ Hence $\lvert v(x) + c \rvert \leq w(x)$ as required. \quad 
\end{proof}

\subsection{Poisson Equation and Invariant \fpmc}
\label{sec:pe}
In this section, we shall prove (see Theorem \ref{th:pe}) that under some sufficient conditions (see Assumption \ref{asmp:suff}) there exist a solution $(\rho,h) \in \mathbb R \times \mathscr B_w$ to the \emph{Poisson Equation} for some fixed $c \in \mathscr B_w$ and a real-valued \ftk{} $\mathcal R$ on $\mathscr B_w$,
\begin{align}
 c + \mathcal R(h) = \rho + h \label{eq:pe}
\end{align}
and \emph{invariant \fpm{}}, $\nu$, satisfying
\begin{align}
 \nu(c + \mathcal R(v)) = \nu(v) + \rho, \forall v \in \mathscr B_w.
\end{align}
As in the theory of MCPs, both Poisson equation and invariant \fpm{} play important roles in studying the stability properties of \ftk{s} and the optimization of the average risk (see Section \ref{sec:ar}).

We apply mainly the same trick used in \cite{hairer2011yet}. Consider an auxiliary weight function $W: \mathbf X \maps [0,\infty)$, which is a real-valued $\mathcal B(\mathbf X)$-measurable function and let $w(x) = 1 + \beta W(x)$ with some positive real number $\beta$.

\begin{assumption} \label{asmp:suff}
Let $\mathcal R$ be real-valued risk map. There exists a function $W: \mathbf X \maps [0,\infty)$ which is $\mathcal
  B(\mathbf X)$-measurable, constants $K \geq 0$, $\gamma \in (0,1)$, 
  $\alpha \in (0,1)$, $\alpha_0 \in (0,\alpha)$ and a real-valued \fpm{} $\nu$ such that (i)
\begin{align}
  \overline{\mathcal R^\sharp}_x(W) \leq \gamma W(x) + K, \forall x \in
  \mathbf X \label{eq:ly}
\end{align}
where $\overline{\mathcal R^\sharp}_x(W)$ is calculated on the space $\mathscr B_{w} = \mathscr B_{1 +\beta W}$, $\beta := \alpha_0/K$, and (ii) 
\begin{align}
 \inf_{x \in B} \left\{ \mathcal R_x(v) - \alpha \nu(v) - \mathcal R_x(u) +
 \alpha \nu(u) \right\} \geq 0 \label{eq:doeb}
\end{align}
whenever $v \geq u \in \mathscr B_{w}$, where $B = \{ x \in \mathbf X: W(x) \leq R \} \in \mathcal
   B(\mathbf X)$ for some $R > 2K/(1-\gamma)$.
\end{assumption}

{\em Remark.} If $\mathcal R$ is coherent, then $\overline{\mathcal R^\sharp} = \mathcal R$ in \eqref{eq:ly}. Specifically, if $\overline{\mathcal R^\sharp} = P$, given that the transition kernel $P$ is considered also as a linear operator (see the remark below Definition \ref{def:rmaps}), then $\nu$ in \eqref{eq:doeb} is a probability measure and \eqref{eq:doeb} is equivalent to $\inf_{x \in B} \left\{ P_x(A) - \alpha \nu(A) \right\} \geq 0, \forall A \in \mathcal B(\mathbf X)$. Therefore, the above assumption is reduced to the condition used in \cite{hairer2011yet}, where it is proven to be in some sense equivalent to the classical geometric ergodicity condition stated in \cite{meyn1993markov}.

In the following, we connect the generalized Doeblin's condition in \eqref{eq:doeb} with the classical one by \emph{subgradients} of risk maps \cite{svindland2009subgradients}. Define the subgradient at state $x \in \mathbf X$ and function $u \in \mathscr B_w$ for a real-valued risk map $\mathcal R$ as follows, 
\begin{align*}
 \delta \mathcal R_x(u):= \left\{ g \left| 
\begin{array}{l}
 g \textrm{ is $\mathcal B(\mathbf X)$-measurable and } \int \lvert g \rvert w dP_x < \infty, \\
 \mathcal R_x(v) \geq \mathcal R_x(u) + \int g(v-u) d P_x, \forall v \geq u \in \mathscr B_w 
\end{array} \right. \right\}
\end{align*}
where $P_x$ denotes the transition probability measure from state $x \in \mathbf X$. 
\begin{proposition}\label{prop:doeb}
Suppose the transition kernel $P$ satisfies that there exists some positive constant $\beta$, set $B\in \mathcal B(\mathbf X)$ and probability measure $\mu$ such that $\inf_{x \in B} \left\{ P_x(A) - \beta \mu(A) \right\} \geq 0, \forall A \in \mathcal B(\mathbf X)$. Assume further that there exists $g(x,u) \in \delta \mathcal R_x(u)$ and positive constant $\epsilon > 0$ such that  $g(x,u) \geq \epsilon$ for all $x \in B$ and $u \in \mathscr B_w$. Then \eqref{eq:doeb} holds for $\nu = \mu$ and $\alpha = \epsilon \beta$. 
\end{proposition}
\begin{proof}
 By definition, we have for each $x \in B$ and $u \in \mathscr B_w$
\begin{align*}
 \mathcal R_x(v) \geq \mathcal R_x(u) + \int g(x,u)(v-u) dP_x \geq R_x(u) + \epsilon \beta \mu(v -u). 
\end{align*}
Then setting $\nu =\mu$ and $\alpha = \epsilon \beta$, \eqref{eq:doeb} holds. \quad 
\end{proof}

{\em Remark.} Subgradients are originally defined for convex risk maps. For concave risk maps $\mathcal R$, however, we can consider its convex counterpart $\tilde{\mathcal R}(v) := - \mathcal R(- v)$. If $\tilde{\mathcal R}$ satisfies the generalized Doeblin's condition with some \fpm{} $\tilde \nu$, then it is easy to see that $\mathcal R$ also satisfies the condition with the \fpm{} $\nu(v) : = - \tilde \nu(-v)$. 


Now we state the contraction theorem under $w$-seminorm.
\begin{theorem}\label{th:semidiscount}
 Suppose Assumption \ref{asmp:suff} holds. Then there exists a $\bar \alpha \in [0,1)$ such that
$\lVert \mathcal R(v) - \mathcal R(u) \rVert_{s,w} \leq \bar \alpha \lVert v - u \rVert_{s,w},$
for all $v$ and $u$ in $\mathscr B_{w}$.
\end{theorem}

{\em Proof}.
Clearly, the assertion is equivalent to $\lVert \mathcal R(v + u) - \mathcal R(u) \rVert_{s,w} \leq \bar \alpha \lVert v \rVert_{s,w}$, $\forall v,u \in \mathscr B_w$. Suppose $\lVert v \rVert_{s,w} \leq C$. Lemma \ref{lm:semi} suggests that we can always find a real value $c$ such that $\lVert v + c \rVert_{w} \leq C$. Since adding any constant to $v$ will not change the values of both sides of the required inequality, without loss of generality, we assume $\lVert v \rVert_{w} \leq C$. Hence, $\lvert v(x) \rvert \leq \lVert v \rVert_{w} w(x) \leq C w(x)$, $\forall x \in \mathbf X$. By definition and Proposition \ref{prop:sharptop}, we have $\forall x \in \mathbf X$,
\begin{align}
 &\lvert \mathcal R_x(v + u) - \mathcal R_x(u) \rvert \leq \overline{\mathcal R^\sharp}_x( \lvert v  \rvert) \leq \lVert v  \rVert_w \overline{\mathcal R^\sharp}_x(w) = C \left( 1 + \beta \overline{\mathcal R^\sharp}_x(W) \right), \label{eq:1}
\end{align}
where the equality is obtained by using Proposition \ref{prop:um}(iii).

We first assume $W(x) + W(y) \geq R$ and set $\gamma_0 := \gamma +\frac{2K}{R} < 1$ and $\gamma_1 := \frac{2 + \beta R \gamma_0}{2+ \beta R} \in (\gamma_0,1)$. 
\eqref{eq:1} yields,
\begin{align}
 & \lvert \mathcal R_x(v + u) - \mathcal R_x(u) -  \mathcal R_y(v + u) + \mathcal R_y(u) \rvert \nonumber\\
\leq & \lvert \mathcal R_x(v + u) - \mathcal R_x(u) \rvert + \lvert \mathcal R_y(v + u) - \mathcal R_y(u) \rvert \nonumber\\
 \leq & C  \left( 2 + \beta \overline{\mathcal R^\sharp}_x(W) + \beta \overline{\mathcal R^\sharp}_y(W)  \right) 
 \leq C \left( 2 + \beta \gamma W(x) + \beta \gamma W(y) + 2 \beta K  \right) \nonumber \\
 \leq& C \left( 2 + \beta \gamma_0 W(x) + \beta \gamma_0 W(y)  \right) \leq  C \left( 2\gamma_1 + \beta \gamma_1 W(x) + \beta \gamma_1 W(y)  \right) \nonumber \\ 
= & C \gamma_1 d_{\beta}(x,y), \label{eq:cond1}
\end{align} 
where the last inequality is due to fact that $\frac{2(1-\gamma_1)}{\beta(\gamma_1-\gamma_0)} = R \leq W(x) + W(y)$. 

Now consider $W(x) + W(y) \leq R$. Thus $x, y \in B$. Define a new \ftk{} $\tilde{\mathcal R}_x(v) := \frac{1}{1-\alpha}\mathcal R_x(v) - \frac{\alpha}{1-\alpha} \nu(v)$.
It is easy to verify that $\tilde{\mathcal R}$ is valid \ftk{} on $B$. In
fact, the monotonicity is guaranteed by Assumption
\ref{asmp:suff}(ii). Hence, by replacing $\mathcal R$ with $\tilde{\mathcal R}$
\eqref{eq:1} holds for all $x,y \in B$, which yields
\begin{align*}
 &\lvert \mathcal R_x(v + u) - \mathcal R_x(u) - \mathcal R_y(v + u) + \mathcal
 R_y(u) \rvert \\
 =& (1 - \alpha) \lvert \tilde{\mathcal R}_x(v + u) -  \tilde{\mathcal
   R}_x(u) -  \tilde{\mathcal R}_y(v + u) + \tilde{ \mathcal
 R}_y(u) \rvert \\
 \leq &  (1 - \alpha) C \left(2 + \beta \overline{\tilde{\mathcal R}^\sharp}_x(W)
 + \beta \overline{\tilde{\mathcal R}^\sharp}_y(W) \right).
\end{align*}
On the other hand $\overline{\tilde{\mathcal R}^\sharp}_x(W) \leq
 (1-\alpha)^{-1}\overline{{\mathcal R}^\sharp}_x(W)$, since $W \geq
0$. Hence,
\begin{align*}
  (1 - \alpha) C (2 + \beta \overline{\tilde{\mathcal R}^\sharp}_x(W)
 + \beta \overline{\tilde{\mathcal R}^\sharp}_y(W) ) \leq &  2 (1 - \alpha) C
 + \beta C \left(2 + \beta \overline{{\mathcal R}^\sharp}_x(W)
 + \beta \overline{{\mathcal R}^\sharp}_y(W) \right) \\
 \leq &  2 (1 - \alpha) C + \beta C (\gamma W(x) + \gamma W(y) + 2 K).
\end{align*}
Since $\beta = \alpha_0 /K$ for some $\alpha_0 \in (0,\alpha)$, setting 
$\gamma_2 := (1 - \alpha + \alpha_0) \vee \gamma \in (0,1)$ yields,
\begin{align}
 & \lvert \mathcal R_x(v + u) - \mathcal R_x(u) - \mathcal R_y(v + u) + \mathcal
 R_y(u) \rvert \nonumber \\
 \leq & 2 C (1 - \alpha + \alpha_0) + C \gamma \beta (W(x) +
 W(y)) \leq C \gamma_2 d_\beta(x,y). \label{eq:cond2} 
\end{align}
Hence, setting $\bar \alpha := \gamma_1 \vee \gamma_2 < 1$, \eqref{eq:cond1} and \eqref{eq:cond2} imply,
$$
 \lvert \mathcal R_x(v + u) - \mathcal R_x(u) - \mathcal R_y(v + u) + \mathcal
 R_y(u) \rvert 
\leq \lVert v \rVert_{s,w} \bar \alpha d_\beta(x,y). \eqno\endproof
$$

\begin{lemma}\label{lm:mu12}
Let $\mathcal T: \mathscr B_w \maps \mathscr B_w$ be an operator satisfying $\mathcal T(v + c) = \mathcal T(v) + c, \forall v \in \mathscr B_w, c \in \mathbb R$ and define its iteration as $\mathcal T^t(v) := \mathcal T(\mathcal T^{t}(v)), t = 2,3,\ldots$ Suppose $\mathcal T$ satisfies furthermore that for all $ v, u \in \mathscr B_w$, (i) $\lVert \mathcal T(v) -\mathcal T(u) \rVert_{s,w} \leq \bar \alpha \lVert v - u \rVert_{s,w}$ for some $\bar \alpha \in [0,1)$; and (ii) $\sup_{t \geq 1} \lVert \mathcal T^t(v) - \mathcal T^t(u) \rVert_w \leq A \lVert v - u \rVert_w$ with some constant $A \in \mathbb R_+$. Then for arbitrary probability measure $\mu_0 \in \mathscr M_w$ and $v,u \in \mathscr B_w$, 
\begin{align}
 \lim_{s \rightarrow \infty} \sup_{t \geq s} \lvert \mu_0 [\mathcal T^{t}(v) - \mathcal T^{t}(u)] - \mu_2[\mathcal T^{s}(v) - \mathcal T^{s}(u)] \rvert = 0. \label{eq:discounttt}
\end{align}
\end{lemma}
\begin{proof}
Let $v_t := \mathcal T^{(t)} (v)$ and $u_t := \mathcal T^{(t)} (u)$, $t = 1, 2, \ldots$ For $v,u \in \mathscr B_w$, without loss of generality, we assume that $\lVert v - u \rVert \leq C$, where $C$ is a positive real number. Hence, for any $t \geq s$, 
\begin{align}
 &\sup_{\lVert v - u \rVert_w \leq C} \lvert \mu_0 [\mathcal T^{t}(v) - \mathcal T^{t}(u)] - \mu_0[\mathcal T^{s}(v) - \mathcal T^{s}(u)] \rvert \nonumber \\
\textrm{(Lemma \ref{lm:semi})} \ = & \sup_{\lVert v -u\rVert_{s,w} \leq C} \lvert \mu_0 [\mathcal T^{t}(v) - \mathcal T^{t}(u)] - \mu_0[\mathcal T^{s}(v) - \mathcal T^{s}(u)] \rvert \nonumber\\
\textrm{(by (i))\ } 
\leq & \sup_{\lVert v_1 - u_1 \rVert_{s,w} \leq \bar \alpha C}  \lvert \mu_0[\mathcal T^{t-1}(v_1) - \mathcal T^{t-1}(u_1)] - \mu_0[\mathcal T^{s-1}(v_1) - \mathcal T^{s-1}(u_1)] \rvert   \nonumber\\
\leq & \sup_{\lVert v_s - u_s \rVert_{s,w} \leq \bar\alpha^s C}  \lvert \mu_0[\mathcal T^{t-s}(v_s) - \mathcal T^{t-s}(u_s)] - \mu_0[v_s - u_s] \rvert \nonumber\\
\textrm{(Lemma \ref{lm:semi})} \ = & \sup_{\lVert v_s - u_s \rVert_{w} \leq \bar \alpha^s C} \lvert \mu_0[\mathcal T^{t-s}(v_s) - \mathcal T^{t-s}(u_s)] - \mu_0[v_s - u_s] \rvert. \nonumber 
\end{align}
Since $\lVert v_s - u_s \rVert_{w} \rightarrow 0$ as $s \rightarrow 0$, and by (ii),
$\sup_{t\geq s} \lVert \mathcal T^{t-s}(v_s) - \mathcal T^{t-s}(u_s) \rVert_w \leq A \lVert v_s - u_s\rVert_w \rightarrow 0,$
the assertion holds. 
\quad \end{proof}

Let $\{ c_t\}$ be a sequence of functions in $\mathscr B_w$ and $\mathcal R$ be a real-valued risk map. Define $\mathcal F_t(v) := c_t + \mathcal R(v)$, $\mathcal F^{(0)}(v) := v$, and $\mathcal F^{(t+1)}(v) := \mathcal F_{t+1} (\mathcal F^{(t)}(v)), t = 0,1,2,\ldots$

\begin{lemma}\label{lm:rhoc}
 Suppose Assumption \ref{asmp:suff} holds. Then for all $v,u \in \mathscr B_w$, (i) $\sup_{t \geq 1} \lVert \mathcal F^{(t)}(v) - \mathcal F^{(t)}(u) \rVert_w < \infty$ and (ii) $\lim_{t \rightarrow \infty} \frac{1}{t} \lVert \mathcal F^{(t)}(v) - \mathcal F^{(t)}(u) \rVert_w = 0$.
\end{lemma}

\begin{proof}
Due to Assumption \ref{asmp:suff}, by \eqref{eq:ly}, we have $\overline{\mathcal R^\sharp}(w) \leq \gamma w + K'$,
where $K' :=  \beta K + 1 - \gamma$. Thus $\lvert \mathcal F^{(1)}(v) - \mathcal F^{(1)}(u) \rvert \leq \overline{\mathcal R^\sharp}(\lvert v - u \rvert) \leq \lVert v - u \rVert_w \overline{\mathcal R^\sharp}(w) \leq \lVert v - u \rVert_w (\gamma w + K').$ 
By induction w.r.t.~$t$, we have for $t = 2, 3, \ldots$
\begin{align}
 & \lvert \mathcal F^{(t)}(v) - \mathcal F^{(t)}(u) \rvert \leq
  \overline{\mathcal R^\sharp}(\lvert \mathcal F^{(t-1)}(v) - \mathcal F^{(t-1)}(u) \rvert)  \nonumber \\
\leq & \lVert v - u \rVert_w \overline{\mathcal R^\sharp}\left(  \gamma^{t-1} w + K'\sum_{i=0}^{t-2} \gamma^i  \right) \leq \lVert v - u \rVert_w \left( \gamma^t w + K'\sum_{i=0}^{t-1} \gamma^i \right)  \label{eq:bound} 
\end{align}
which implies $\sup_{t \geq 1} \lVert \mathcal F^{(t)}(v) - \mathcal F^{(t)}(u) \rVert_w < \infty$. (ii) is an immediate result of (i). 
\end{proof}

\begin{theorem} \label{th:pe}
 Suppose Assumption \ref{asmp:suff} holds. Then, for each $c \in \mathscr B_w$, (i) the Poisson equation
\eqref{eq:pe} has one solution $(\rho,h) \in \mathbb R \times \mathscr B_w$, where $\rho$ is unique, and (ii) there exists a real-valued \fpm{} $\nu$ such that $\nu(c + \mathcal R(v)) = \nu(v) + \rho$, $\forall v \in \mathscr B_w$.
\end{theorem}
\begin{proof}
(i) Define $\mathcal T_c(\cdot) := c + \mathcal R(\cdot)$. Then by Assumption \ref{asmp:suff} and Theorem \ref{th:semidiscount}, the map $\mathcal T_c: \mathscr B_w \rightarrow \mathscr B_w$ is a contraction under $w$-seminorm, i.e., $\lVert \mathcal T_c (v) - \mathcal T_c (u) \rVert_{s,w} \leq \bar \alpha \lVert v - u \rVert_{s,w}$, for some $\bar \alpha \in [0,1)$. In the following, we extend the fixed-point theorem w.r.t.\ span-seminorm (cf.\ p.\ 321 \cite{arapostathis1993discrete} for bounded $w$) to $w$-seminorm. Let $\widetilde{\mathscr B}_w = \mathscr B_w / \sim$ be the quotient space, which is induced by the equivalence relation $\sim$ on $\mathscr B_w$ defined by $v \sim u$ if and only if there exists some constant $C \in \mathbb R$ such that $v(x) - u(x) = C$ for all $x \in \mathbf X$, endowed with the quotient norm induced by the $w$-seminorm. For $v \in \mathscr B_w$, let $\tilde v$ be the corresponding equivalent class in $\tilde{\mathscr B_w}$ and $\widetilde{\mathcal T_c}: \widetilde{\mathscr B}_w \rightarrow \widetilde{\mathscr B}_w$ be the canonically induced map, i.e., $\widetilde{\mathcal T_c}(\tilde v) := \widetilde{\mathcal T_c(v)}$, $v \in \mathscr B_w$. Since $\mathcal T_c$ is a contraction w.r.t.\  $w$-seminorm on $\mathscr B_w$, $\widetilde{\mathcal T_c}$ is a contraction on $\widetilde{\mathscr B}_w$ and therefore has a unique fixed point. Conversely, it follows that the map $\mathcal T_c$ has a $w$-seminorm fixed point. In other words, there exists $h \in \mathscr B_w$ such that $\lVert \mathcal T_c(h) - h \rVert_{s,w} = 0$ and $\rho_c := \mathcal T_c(h) - h$ is a constant.

Next we show that such $\rho_c$ is unique. Define $\mathcal T_c^t(\cdot) := \mathcal T_c(\mathcal T_c^{t-1}(\cdot)), t = 2,3,\ldots$ Suppose there exits another solution $(\rho',h') \in \mathbb R \times \mathscr B_w$ to the Poisson equation. Then $\mathcal T_c^t(h') = t \rho + h'$ and $\mathcal T_c^t(h) = t \rho_c + h$. However, by Lemma \ref{lm:rhoc} (ii), $\frac{1}{t}\lVert \mathcal T_c^t(h') - \mathcal T_c^t(h) \rVert_w \rightarrow 0$,  
which implies that $\lim_{t \rightarrow \infty} \frac{1}{t}\lVert  t \rho' + h' - t \rho_c - h \rVert_w = 0$. Hence, $\rho' = \rho_c$.

(ii) Let $h \in \mathscr B_w$ be a solution to the Poisson equation and $\mu_0 \in \mathscr M_w$ is a probability measure. Define $\mathcal D_c(\cdot) := \mathcal T_c(\cdot) - \rho_c$ and $\mu_t(\cdot) := \mu_0(\mathcal D_c^t(\cdot))$. Then we have $\forall v \in \mathscr B_w$, 
\begin{align*}
 \mathcal T^{t}_c(v) - \mathcal T^{t}_c(h) - (\mathcal T^{s}_c(v) - \mathcal T^{s}_c(h)) = \mathcal T^{t}_c(v) -  \mathcal T^{s}_c(v) - (t-s) \rho_c = \mathcal D_c^t(v) - \mathcal D_c^s(v).
\end{align*}
Due to \eqref{eq:bound}, the condition (ii) in Lemma \ref{lm:mu12} holds. Hence, setting $u = h$ in Lemma \ref{lm:mu12}, we obtain that 
\begin{align*}
 \lim_{s \rightarrow \infty} \sup_{t \geq s} \lvert \mu_0 [\mathcal T^{t}_c(v) - \mathcal T^{t}_c(h)] - \mu_2[\mathcal T^{s}_c(v) - \mathcal T^{s}_c(h)] \rvert = 0,
\end{align*}
which implies for each $v \in \mathscr B_w$,
\begin{align*}
 \lim_{s \rightarrow \infty} \sup_{t \geq s} \mu_0 [  \mathcal D_c^t(v) - \mathcal D_c^s(v) ] =\lim_{s \rightarrow \infty} \sup_{t \geq s} [ \mu_t (v) - \mu_s(v)] = 0.
\end{align*}
Hence $\mu_t$ converges to a mapping $\mu_\infty: \mathscr B_w \rightarrow \mathbb R$ satisfying $\mu_\infty(\mathcal D_c(v)) = \mu_\infty(v)$, $\forall v \in \mathscr B_w$. Other other hand, it is easy to see that for each $t$, $\mu_t$ is a real-valued \fpm{} except the axiom of centralization. Hence, $\mu_\infty(c + \mathcal R(v)) = \mu_\infty(v) + \rho_c$. Finally, by setting $\nu(v) := \mu_\infty(v) - \mu_\infty(0)$, we obtain the required \fpm{}.
\quad \end{proof}

\section{Applying Risk Measures to MCPs}
\label{sec:armtmcp}

\subsection{\ftkc{s} for MCPs}
We define the \ftk{s} with controls as follows
\begin{definition}\label{def:gen}
  $\mathcal R(v|x,a)$ (simply written as $\mathcal R$) is said to be a \textbf{\ftk} on an MCP $(\mathbf X, \mathbf A,\{\mathbf A(x)|x \in \mathbf X \},Q)$, if (i) for each $(x,a) \in \mathbf K$, $\mathcal R(\cdot|x,a): \mathscr B_w \maps \mathbb R$ is a real-valued \fpm{}; and (ii) for each $v \in \mathscr B_w$, $\mathcal R(v|\cdot)$ is a real-valued $\mathcal B(\mathbf K)$-measurable function.
Furthermore, we define for any $\pi \in \Delta$
\begin{align}
 \mathcal R^\pi(v|x) := \int_{\mathbf A(x)} \pi(\diff a | x)\mathcal R(v|x,a). \label{eq:linear}
\end{align}
\end{definition}
For convenience, we sometimes write $\mathcal R_{x,a}(v) := \mathcal R(v|x,a)$ and $\mathcal R_{x}^\pi(v) := \mathcal R^\pi(v|x)$.

{\em Remark.} Note that since \ftk{s} are subjective representations of objectives transition probabilities, as in the above definition, $\mathcal R$ depends always on the transition model $Q$ of the underlying MCP. It is obvious that the transition kernel $P^\pi$ defined in \eqref{eq:cpi} is a valid \ftk{}. Thus, the concept of a \ftk{} is a generalization of the conditional expectation.  \eqref{eq:linear} in fact assumes that $\mathcal R^\pi$ is linear to the policy $\pi$, which simplifies the optimization problem and is one of the conditions that guarantee the existence of one optimal deterministic policy (``optimal selector'').

In the mathematical finance literature, there exist various ways to extend the CRM to a temporal structure (see e.g.~\cite{
   foellmer2006convex,cheridito2011composition,
   ruszczynski2010risk} and references therein). The definition is usually selected based on applications. 
To compare their subtle differences are out of the scope of this paper. The \ftk{s} defined here are similar to the \emph{risk measure generators} in \cite{cheridito2011composition} and are implicitly Markovian and time-homogeneous (see also \cite{ruszczynski2010risk}), since $\mathcal R$ defined above depends merely on the most recent state and action but not the whole history. 
The risk maps used in this paper are assumed to be Markovian, since in the MCP-framework the underlying stochastic process is Markov, while the assumption of time-homogeneity is due to the fact that since we consider the infinite-horizon criteria (see \eqref{eq:JalLim} and \eqref{eq:avrk}), as in the literature of MCPs, stationary optimal policies are expected. 
Hence, to comply with the MCP-framework, it is sufficient to construct an operator which replaces the conditional expectation determined by the transition model $Q$ and policy $\pi$.

\subsection{Risk Preference}
\label{sec:pre}
It is important to know how to judge the risk-preference of one specific \ftk{} according to the properties of the map. In the following, we introduce the correspondence between the risk-preference and the convexity (concavity) by an intuitive example. Recall that in economics risk represents uncertainties \cite{chavas2004risk}. Thus there is no risk if there is no uncertainty. Suppose we are given two options, the first one is $(v,p;u,1-p)$, $0 < p < 1$, to obtain outcome (reward or cost) $v$ with probability $p$ and outcome $u$ with probability $1-p$; the second one is to obtain outcome $p\cdot v + (1-p) \cdot u$ with probability 1. The first one is riskier than the second one in the sense that though both options have the same mean outcome, the first one includes uncertainty $p$ (and $1-p$) whereas the outcome is certain in the second option. Thus, if the first one is preferred, the agent is risk-seeking; if the second one is preferred, the agent is risk-averse.

In our framework, the preference is determined by the value $\mathcal R_{x,a}(v)$. Note that we want to minimize risk. If the $\mathcal R$ is concave w.r.t.~$v$ at some $(x,a) \in \mathbf K$, 
\begin{align*}
 p \mathcal R_{x,a}(v) + (1-p) \mathcal R_{x,a}(u) \leq \mathcal R_{x,a}(p v + (1-p) u), \forall p \in [0,1], v, u \in \mathscr B
\end{align*}
the riskier option is preferred. Thus, $\mathcal R$ is intuitively risk-seeking at $(x,a)$. Conversely, if $\mathcal R$ is convex w.r.t.~$v$ at $(x,a)$, the agent is risk-averse at $(x,a)$.  

{\em Remark.} The categorization depends on the objective. In the \emph{expected utility theory} \cite{gollier2004economics}, the objective is to maximize utilities. Therefore, the categorization is opposite: concavity means risk-averse and convexity suggests risk-seeking. Several existing risk measures (see Section \ref{sec:ex}) in the literature confirm also the above defined categorizations from intuitive sense. 

In some cases, it is reasonable to take different risk-preferences at different situations, e.g.~at ``safe'' states the bolder policies are taken to obtain high uncertain reward, while at ``dangerous'' states the conservative policies are applied to avoid high uncertain cost. Even when risk-averse policies are applied everywhere, the degree of risk-averse at different states can be also different. Within our framework, mixed risk-preferences can be easily modeled due to the following lemma.
\begin{lemma} Let $\mathcal R, \mathcal R'$ be two risk maps. For any $B \in \mathcal B(\mathbf K)$, we define
$$\tilde{\mathcal R}(v|x,a) := \mathbf{1}_{B}(x,a) \mathcal R(v|x,a) + \mathbf{1}_{B^\textrm{C}}(x,a) \mathcal R'(v|x,a), $$
where $\mathbf{1}_{B}(\cdot)$ denotes the indicator function and $B^\textrm{C}$ is the complementary set of $B$. Then $\tilde{\mathcal R}$ is also a risk map.
\end{lemma}
\begin{proof}
Obviously for each $(x,a) \in \mathbf K$, $\tilde{\mathcal R}(\cdot|x,a)$ is a risk measure, since both $\mathcal R(\cdot |x,a)$ and $\mathcal R'(\cdot|x,a)$ are risk measures. Given $v \in \mathscr B_w$, since $\mathcal R(v |\cdot)$ and $\mathcal R'(v|\cdot)$ are both real-valued $\mathcal B(\mathbf K)$-measurable functions, it follows that $\mathcal R(v |\cdot)$ is also a real-valued $\mathcal B(\mathbf K)$-measurable function. 
\quad \end{proof}

For instance, suppose $\mathcal R$ and $\mathcal R'$ are maps of everywhere risk-averse and risk-seeking respectively. Then, by above lemma, $\tilde{\mathcal R}$ would be the risk map that is risk-averse at $B$ but risk-seeking at $B^C$ (e.g.~$B := \left\{ k \in \mathbf K| c(k)\geq 0\right\}$). Therefore, the lemma gives us the power to model mixed risk-preferences by constructing new maps from existing maps whose properties have been well investigated.

\subsection{Examples}
\label{sec:ex}
There exist several important \ftk{s} in the literature of economics, mathematical finance and control theory. Most of them can be adapted to the framework we introduced above. Since most of the literature consider merely bounded space $\mathscr B$ (or essentially bounded space) except the classical MCPs, in this section, we restrict ourselves to $\mathscr B$ when introducing definitions of different maps. Note that we focus mostly on the literature of optimal control and operation research. For more examples in mathematical finance, we refer to Chapter 4 of \cite{follmer2004stochastic}.
\paragraph{Classical MCPs} 
$\mathcal R_{x,a}(v) := \mathbb E_{x,a}^Q\left[ v \right] = \int_{\mathbf X} Q(\diff y|x,a) v(y)$. 
It is easy to verify that it is coherent and linear to $v$ (therefore risk-neutral).
\paragraph{Entropic map} The name is taken from the literature of CRM \cite{follmer2002convex}. It is intensely also researched in the field of optimal control \cite{chung1987discounted,
  avila1998controlled, borkar2002risk, di2008infinite,coraluppi2000mixed,fleming1997risk,hernandez1996risk,
  marcus1997risk, cavazos2010optimality}. 
\begin{align}
 \mathcal R_{x,a}(v) := \dfrac{1}{\lambda}\log \mathbb E_{x,a}^Q \left[  \exp(\lambda v ) \right] = \dfrac{1}{\lambda} \log \left\lbrace \int_{\mathbf X} Q(\diff y|x,a) \exp( \lambda v(y)) \right\rbrace \label{eq:entropic}
\end{align}
where the risk-sensitive parameter $\lambda \in \mathbb R$ controls the risk-preference of $\mathcal R$: if $\lambda > 0$, $\mathcal R$ is everywhere convex and therefore everywhere risk-averse; if $\lambda < 0$, $\mathcal R$ is everywhere concave and therefore everywhere risk-seeking. It can be also shown that
\begin{align*}
 \lim_{\lambda \rightarrow 0} \dfrac{1}{\lambda} \log \left\lbrace \int_{y} Q(\diff y|x,a) \exp( \lambda v(y)) \right\rbrace = \int_{y} Q(\diff y|x,a) v(y)
\end{align*}
which is equivalent to the classical MCPs. Besides, it has connection to the mean-variance tradeoff scheme via the Taylor expansion at $\lambda = 0$, $
 \frac{1}{\lambda} \log \mathbb E \exp(\lambda Z ) = \mathbb E Z + \lambda \textrm{Var}(Z) + O(\lambda^2)$
where $Z$ denotes an arbitrary bounded random variable. Suppose that risk is measured by variance. The objective is to minimize risk $\mathcal R^\pi$. Therefore, if $\lambda > 0$, the variance is avoided, the agent is risk-averse. On the contrary, if $\lambda < 0$, the variance is preferred, the agent is intuitively risk-seeking. These intuitions coincide the categorization based on the convexity (concavity) of $\mathcal R$. 

\paragraph{Robust risk maps} Iyengar (2005) \cite{iyengar2005robust} introduced the framework of \emph{robust dynamic programming}, by which he argues that in some applications the transition model $Q$ cannot be inferred exactly. Instead, he employs a set of transition probabilities, $\mathcal P$, which contains all possible ``ambiguous'' transition kernels. In order to gain the ``robustness'', the worst cost is considered, adapted 
in our framework, 
\begin{align}
\mathcal R_{x,a}(v) := \sup_{Q(\cdot|x,a)\in \mathcal P_{x,a}} \mathbb E_{x,a}^Q v = \sup_{Q(\cdot|x,a) \in \mathcal P_{x,a}} \int_{y \in \mathbf X} Q(\diff y|x,a) v(y). \label{eq:rob}
\end{align}
 It is apparent that $\mathcal R$ is coherent. We can verify that $\mathcal R$ is everywhere convex and therefore risk-averse, which coincides the intuition that the worst scenario is considered. One special case of the robust dynamic programming was the \emph{minimax control} (see e.g.~\cite{coraluppi2000mixed}), which also considers the worst scenario with finite state space, $\mathcal R_{x,a}(v) := \max_{Q(y|x,a) > 0} v(y).$ It is also notable that each coherent risk map has one dual presentation of the form \eqref{eq:rob} under some regularity conditions for the set $\mathcal P_{x,a}$ (see e.g.~\cite{Delbaen_2000} for essentially bounded spaces and \cite{svindland2009convex} for unbounded ones).

\paragraph{Mean-semideviation trade-off} \cite{ogryczak1999stochastic,ruszczynski2006optimization} This risk map considers only the trade-off between the one-step conditional mean and semideviation rather than the deviation of the whole Markov chain \cite{sobel1982variance, filar1989variance}.
\begin{align}
 \mathcal R_{x,a}(v) := \mathbb E_{x,a}^Q[v] + \lambda \left[ \mathbb E_{x,a}^Q (v - \mathbb E_{x,a}^Q[v])^r_+ \right]^{1/r} \label{eq:meansemi}
\end{align}
where $r \geq 1$ and $\lambda \in \mathbb R$ denotes the risk-preference parameter which controls the risk preference of $R$: if $\lambda >0$, $R$ is risk-averse; if $\lambda < 0$, $R$ is risk-seeking. Setting $r=2$, this map can be viewed as an approximation of the mean-variance tradeoff scheme defined in \cite{filar1989variance}.

\paragraph{Prospect theory and Choquet integral} The idea is based on \emph{non-additive measures} (Denneberg, 1994, \cite{denneberg1994non}) and \emph{capacities} (Choquet, 1953, \cite{choquet1953theory}) that generalize standard probability measures (p.m.). \emph{Non-additive p.m.} (for definition see \cite{denneberg1994non}) can be analogously extended to fit the MCP framework with \emph{conditional non-additive p.m.}~$\Phi$, which satisfies that (i) for each $(x,a) \in \mathbf K$, $\Phi(\cdot|x,a)$ is a non-additive p.m.~and (ii) for each $B \in \mathcal B(\mathbf X)$, $\Phi(B|\cdot)$ is $\mathcal B(\mathbf K)$-measurable. Then we define the conditional \emph{Choquet integral} as 
\begin{align*}
 \int_{\mathbf X}^{Ch} v(y) \Phi(\diff y|x,a) := \int_{-\infty}^0 [\Phi(v(y)>t|x,a) - 1] \diff t +  \int_0^\infty \Phi(v(y)>t|x,a) \diff t.
\end{align*}
It is easy to see that for each conditional non-additive p.m.~$\Phi$, Choquet integral is a homogeneous risk map but not necessarily convex everywhere. The well-known \emph{prospect theory} \cite{tversky1992advances} in behavioral economics, which is applied to interpret human behaviors of mixed risk-preferences, can be in fact represented as a Choquet integral \cite{kahneman1979prospect, tversky1992advances}. In a recent study, Chateauneuf and Cohen (2008) \cite{chateauneuf2008cardinal} apply the Choquet integral to model the \emph{subjective expected utility} \cite{savage1972foundations}.

\section{Optimal Risk-sensitive Control}
\label{sec:optrisk}
Suppose $\mathcal R$ is a \ftk{} and $\bs \pi \in \Pi_M$. Analogous to the classical MCPs, we consider the following three objectives: (a) \textbf{$T$-stage total risk},
\begin{align}
 J_T(x,\bs \pi) := c^{\pi_0}(x) + \mathcal R^{\pi_0}_{x} \left( c^{\pi_1} + \mathcal R^{\pi_1} \left(c^{\pi_2} + \ldots + \mathcal R^{\pi_{T-1}}(c^{\pi_T}) \ldots \right) \right) \label{eq:JT}
\end{align}
(b) the \textbf{discounted total risk} as
\begin{align}
 J_{\alpha}(x,\bs \pi) := \lim_{T \rightarrow \infty} J_{\alpha,T}(x,\bs \pi) \label{eq:JalLim}
\end{align}
where the discounted $T$-stage risk is defined as follows,
\begin{align*}
 J_{\alpha,T}(x,\bs \pi) : = c^{\pi_0}(x) + \alpha \mathcal R^{\pi_0}_{ x} \left( c^{\pi_1} + \alpha \mathcal R^{\pi_1} \left(c^{\pi_2} + \ldots + \alpha \mathcal R^{\pi_{T-1}}(c^{\pi_T}) \ldots \right) \right) 
\end{align*}
and (c) the \textbf{average risk}
\begin{align} \label{eq:avrk}
 J(x, \bs \pi) := \limsup_{T \rightarrow \infty} \dfrac{1}{T}J_T(x, \bs \pi), \bs \pi \in \Pi_M, x \in \mathbf{X}.
\end{align}
The optimal control problems for above three risk-sensitive objectives are to minimize the risk among all Markov policies 
\begin{align*}
 J_T^*(x) := \inf_{\bs \pi \in \Pi_M} J_T(x,\bs \pi),\quad J_\alpha^*(x) := \inf_{\bs \pi \in \Pi_M} J_{\alpha}(x,\bs \pi), \textrm{and} \quad  
J^*(x) := \inf_{\bs \pi \in \Pi_M} J(x,\bs \pi).
\end{align*}
In the rest of this section, for convenience, all the three objectives can be considered as functions on $\mathbf X$ within the space $\mathscr B_w$ by using the notations  $J_T(\bs \pi), J_\alpha(\bs \pi)$ and $J(\bs \pi)$, as well as $J^*_T$, $J^*_\alpha$ and $J^*$.

The finite-stage total risk can be solved by \emph{dynamic programming} (cf.~\cite{ruszczynski2010risk}). 
In this paper, we focus on the discounted total risk and average risk, which will be discussed in detail in
Section \ref{sec:disc} and \ref{sec:ar} respectively.

In the rest of this section, we shall frequently use the following operators
\begin{align}
 \mathcal F^\pi_\alpha(v) := c^\pi + \alpha \mathcal R^\pi(v), \quad \mathcal F_\alpha(v) := \inf_{\bs \pi \in \Pi_M} \mathcal F^\pi_\alpha(v), \quad v \in \mathscr B_w  \label{eq:F}
\end{align}
where $\alpha \in [0,1]$. If $\alpha = 1$, we simply write them as $\mathcal
F^\pi$ and $\mathcal F$ accordingly. The operators $(\mathcal F^\pi_\alpha)^t$, $t \in \mathbb N$, is defined iteratively as $(\mathcal F^\pi_\alpha)^0(v) := v$, and $(\mathcal F^\pi_\alpha)^t(v) := \mathcal F^\pi_\alpha((\mathcal F^\pi_\alpha)^{t-1}(v))$, $t = 1,2,\ldots$, while $\mathcal F^t_\alpha$ is defined analogously.

Analogous to classical MCPs, we need further assumptions to guarantee the existence of the ``selector'' in the optimization problem.  
\begin{assumption} \label{asmp:sub}
For each $x \in \mathbf X$,
 \begin{enumerate}[(a)]
  \item the cost function $c(x,a)$ is lower semi-continuous on $\mathbf A(x)$,
  \item the action space $\mathbf A(x)$ is compact, and
  \item the function $u'(x,a) := \mathcal R_{x,a}(u)$ is continuous in $a \in \mathbf A(x)$ for any $u \in \mathscr B_w$.
 \end{enumerate}
\end{assumption}

\begin{proposition}
\label{prop:sub}
 Suppose $\mathcal R$ is a \ftk{} satisfying Assumption \ref{asmp:sub}. Then, for all $v \in \mathscr B_w$ and $x \in \mathbf X$, there exists a deterministic policy $f \in \Delta_D$, such that for any $\alpha \in [0,1]$, $c^f(x) + \alpha \mathcal R^f(v|x) = \mathcal F_\alpha(v|x) = \inf_{\pi \in \Delta} \left\{ c^\pi(x) + \alpha \mathcal R^\pi(v|x) \right\}.$
\end{proposition}
\begin{proof}
Apparently, $\mathcal F_\alpha(v|x) = \inf_{a \in \mathbf A(x)} \left\{ c(x,a) + \alpha \mathcal R(v|x,a) \right\}$ for each $x \in \mathbf X$. By Assumption \ref{asmp:sub}(a) and (c), the function
$u(x,a) := c(x,a) + \alpha \mathcal R(v|x,a), \alpha \in [0,1]$, 
is lower semi-continuous in $a \in \mathbf A(x)$ for each $x \in \mathbf X$. Hence, by Assumption \ref{asmp:sub}(b) and Lemma 8.3.8(a) \cite{hernandez1999further}, the optimal selector exists.
\quad \end{proof}

\subsection{Discounted Total Risk}
\label{sec:disc} 
In economics, the time-discount is added to reflect the ``time-value'' of
outcomes: the outcome to be gained in the future is less valuable than the
same amount of outcome obtained now. It has
similar effects when cost is concerned. Due to its good mathematical
properties, exponential discounting scheme, where the cost $c_t$ is
multiplied with the time-discount $\alpha^t$, is widely applied in economics, finance as well
as in MCPs (see $S_\alpha$ in
\eqref{eq:totalreward}). In fact, the natural extension of classical discounted MCPs is as follows,
\begin{align*}
 D_\alpha(\bs \pi) = c^{\pi_0} + \mathcal R^{\pi_0} \left( \alpha c^{\pi_1} + \mathcal R^{\pi_1}\left( \alpha^2 c^{\pi_2} + \ldots + \mathcal R^{\pi_{T-1}}(\alpha^T c^{\pi_T} + \ldots ) \ldots \right) \right). 
\end{align*}
However, since the \ftk{} $\mathcal R$ is not necessarily homogeneous, there
might not exist a stationary policy when optimizing $D_\alpha$. Indeed,
Chung \& Sobel (1987) \cite{chung1987discounted} proved that for the entropic
risk maps, which are not homogeneous, the optimal policy might not be
stationary if
$D_\alpha(\bs \pi)$ is optimized w.r.t.~$\bs \pi$, though $D_\alpha$ is
well-defined if $\alpha \in [0,1)$.

In our definition, discount factor $\alpha$ is multiplied with $\mathcal
R$, which has the same ``time-discount'' effect, where the risk rather than
the immediate cost is discounted. Moreover, it is easy to see that, if $\mathcal R$ is homogeneous,
$D_\alpha$ is equivalent to $J_\alpha$, the discounted total risk under our
definition. Therefore, $D_\alpha$ defined for any homogeneous risk map is
merely a special case of our definition. Specifically, the classical
discounted MCP is indeed a special case of our defined discounted total risk, since
it is homogeneous. Ruszczy{\'n}ski (2010) \cite{ruszczynski2010risk} used
$D_\alpha$ as the objective function, which was solved by a value iteration
algorithm, since merely the coherent risk maps (both homogeneous and convex)
were considered. Besides, in the proof of the value iteration algorithm, he
used the representation theorem,  which
is valid merely for coherent risk maps. On the contrary, we will see later
that the objective $J_\alpha$ allows a value iteration algorithm for general
risk maps. Therefore, we apply $J_\alpha$ rather than $D_\alpha$.

Let $\bs \pi = [\pi_0, \pi_1, \ldots] \in \Pi_M$ be one Markov policy.  Using the operator $\mathcal F^\pi$ defined in \eqref{eq:F}, we have
$$J_{\alpha,T}(\bs \pi) = \mathcal F^{\pi_0}_\alpha (\mathcal F^{\pi_1}_\alpha \ldots \mathcal F^{\pi_{T-1}}_\alpha(c^{\pi_T}) \ldots) = \mathcal F^{\pi_0}_\alpha (\mathcal F^{\pi_1}_\alpha \ldots \mathcal F^{\pi_{T-1}}_\alpha(\mathcal F^{\pi_{T}}_\alpha(0)) \ldots).$$ We first show that under some assumption, the limit in \eqref{eq:JalLim} exists. 
\begin{assumption}\label{assmp:disc}
 There exist nonnegative constants $\overline c$ and $\overline w$, with $1 \leq \overline w < 1/\alpha$ and a weight function $w \geq 1$ such that $\forall x \in \mathbf X$:(a) $\sup_{a \in \mathbf A(x)} \lVert c(x,a) \rVert \leq \overline c w(x) $ and (b) $\sup_{a \in \mathbf A(x)} \overline{\mathcal R^\sharp}(w|x,a) \leq \overline w w(x) $, where $\overline{\mathcal R^\sharp}$ is defined in Section \ref{sec:mod}.
\end{assumption}

Applying Proposition \ref{prop:wbar} we obtain the following proposition.
\begin{proposition}\label{prop:disc2}
 Suppose Assumption \ref{assmp:disc}(b) holds. Then for all $\pi \in \Delta, v , u \in \mathscr B_w$, $  \lVert \mathcal F^\pi_\alpha(v) - \mathcal F^\pi_\alpha(u) \rVert_w \leq \overline w \alpha \lVert v - u \rVert_w$.
\end{proposition}

\begin{lemma}\label{lm:disc2}
 Suppose Assumption \ref{assmp:disc}(a) and (b) holds. Then (i)$\lim_{T\rightarrow \infty} J_{\alpha,T}(\bs \pi)$ exists in $\mathscr B_w$, for all $\bs \pi \in \Pi_M$, and (ii) $J_{\alpha,T}(\bs \pi) = \lim_{T \rightarrow \infty} \mathcal F^{\pi_0}_\alpha (\mathcal F^{\pi_1}_\alpha \ldots \mathcal F^{\pi_{T}}_\alpha(v) \ldots)$, $\forall v \in \mathscr B_w$.
\end{lemma}
\begin{proof}
 By Assumption \ref{assmp:disc}(a), $ c^\pi \in \mathscr B_w$ for all $\pi \in \Delta$. Iterating Proposition \ref{prop:disc2},
\begin{align}
 & \lVert J_{\alpha,T+1}(\bs \pi) - J_{\alpha,T}(\bs \pi) \lVert_w \nonumber \\
= &  \lVert \mathcal F^{\pi_0}_\alpha (\mathcal F^{\pi_1}_\alpha \ldots \mathcal F^{\pi_{T}}_\alpha(c^{\pi_{T+1}}) \ldots) - \mathcal F^{\pi_0}_\alpha (\mathcal F^{\pi_1}_\alpha \ldots \mathcal F^{\pi_{T-1}}_\alpha(\mathcal F^{\pi_{T}}_\alpha(0)) \ldots) \lVert_w \label{eq:tmp1}\\
 \leq &(\overline w \alpha)^T \lVert c^{\pi_{T+1}} \rVert_w. \nonumber 
\end{align}
By Assumption \ref{assmp:disc}(a), $0 \leq \overline w \alpha< 1$. Thus $J_{\alpha,T+1}(\bs \pi) \rightarrow J_{\alpha,T}(\bs \pi)$ and $\lVert J_{\alpha,T+1}(\bs \pi) \lVert_w < \dfrac{\overline c}{1 - \overline w \alpha}$ is bounded and thus in $\mathscr B_w$. 
(ii) is straightforward by replacing $c^{\pi_{T+1}}$ with $v$ in \eqref{eq:tmp1}. 
\quad \end{proof}

\begin{lemma}
 Suppose Assumption \ref{assmp:disc} and \ref{asmp:sub} hold. Then $\lVert\mathcal  F_\alpha(v) - \mathcal F_\alpha(u) \rVert_w \leq \overline w \alpha \lVert v - u \rVert_w, 0 \leq \bar w \alpha < 1.$
\end{lemma}
\begin{proof}
 By Proposition \ref{prop:sub}, the optimal selector $f^*$ always exists for all $v \in \mathscr B_w$. Let $f_v$ be the optimal selector for $v$ and $f_u$ be the optimal selector for $u$. Thus
\begin{align*}
  \mathcal F_\alpha(v) - \mathcal F_\alpha(u) \leq \mathcal F^{f_u}_\alpha(v) - \mathcal F_\alpha^{f_u}(u) \leq \overline w \alpha \lVert v - u \rVert_w w
\end{align*}
where the latter inequality is due to Proposition \ref{prop:disc2}. On the other hand,
\begin{align*}
  \mathcal F_\alpha(u) - \mathcal F_\alpha(v) \leq \mathcal F^{f_v}_\alpha(u) - \mathcal F_\alpha^{f_v}(v) \leq \overline w \alpha \lVert v - u \rVert_w w.
\end{align*}
Hence, $\lVert \mathcal F_\alpha(v) - \mathcal F_\alpha(u) \rVert_w \leq \overline w \alpha \lVert v - u \rVert_w$.
\quad \end{proof}

Since $\mathcal F_\alpha$ is a contracting map and $\mathscr B_w$ is complete, by the Banach fixed point theorem, it has the unique fixed point in $\mathscr B_w$ such that
\begin{align}
 v^*(x) =\mathcal F_\alpha(v^*|x) = \min_{a \in \mathbf A(x)} \left\{ c(x,a) + \alpha \mathcal R(v^*|x,a)\right\}.
\end{align}
The final task is to show that $v^* = J^*_\alpha$, the optimal discounted total risk.

\begin{theorem}\label{th:disc}
 Suppose Assumption \ref{assmp:disc} and \ref{asmp:sub} hold. Then $v^* = J^*_\alpha$.
\end{theorem}
\begin{proof}
 First, we show $v \geq \mathcal F_\alpha v$ implies $v \geq J^*_\alpha$. By Proposition \ref{prop:sub}, we assume that the optimal selector for $v$ is $f$. Hence, we obtain
\begin{align*}
 v \geq \mathcal F^f_\alpha v \geq \mathcal F^f_\alpha \mathcal F^f_\alpha v \geq \mathcal F^f_\alpha \mathcal F^f_\alpha \mathcal F^f_\alpha(v) \geq  (\mathcal F^f_\alpha)^\infty (v) = J_\alpha(f^\infty) \geq J^*_\alpha
\end{align*}
where the equality is due to Lemma \ref{lm:disc2}(ii). Second, we show that $v \leq \mathcal F_\alpha v$ implies $v \leq J^*_\alpha$. Indeed, let $\bs \pi = [\pi_0,\pi_1, \ldots ] \in \Pi_M$ be an arbitrary Markov random policy. Then, for all $\pi \in \Delta$, $v \leq \mathcal F_\alpha v \leq \mathcal F_\alpha^\pi v$. Hence,
\begin{align*}
 v \leq \mathcal F_\alpha^{\pi_0}(v) \leq \mathcal F_\alpha^{\pi_0}\mathcal F_\alpha^{\pi_1}(v) \leq \mathcal F_\alpha^{\pi_0}\mathcal F_\alpha^{\pi_1}\mathcal F_\alpha^{\pi_2}(v) \leq \mathcal F_\alpha^{\pi_0}\ldots \mathcal F_\alpha^{\pi_T}(v) \uparrow J_\alpha(\bs \pi).
\end{align*}
The limit is due to Lemma \ref{lm:disc2}(ii). Since $\bs \pi$ can be arbitrarily chosen, it follows
$ v \leq \inf_{\bs \pi \in \Pi_M} J_\alpha(\bs \pi) = J_\alpha^*.$ Combining above two steps yields $v^* = J^*_\alpha$. 
\quad \end{proof}

\begin{corollary}
 Under Assumption \ref{assmp:disc} and \ref{asmp:sub}, there exists a stationary deterministic policy $f^* \in \Delta_D$ such that 
$J_\alpha^* = J_\alpha((f^*)^\infty)$.
\end{corollary}

\subsection{Average Risk}
\label{sec:ar}

In the first part of this section, we extend the Assumption \ref{asmp:suff} to the MCP-framework and compare it with assumptions in the existing literature of classical MCPs and entropic maps. In the second part, we prove the existence of solutions to the \emph{average risk optimality equation}. 

\begin{assumption}\label{assp:semi}
 There exists a function $W: \mathbf X \maps [0,\infty)$ which is $\mathcal
  B(\mathbf X)$-measurable and constants $K\geq0$, $\gamma \in (0,1)$, $\alpha \in (0,1)$ and a \fpm{} $\nu$ s.t.~(i)
\begin{align}
  \overline{{\mathcal R}^\sharp}_{x,a}(W) \leq \gamma W(x) + K, \forall (x,a) \in
  \mathbf K \label{eq:ly2}
\end{align}
where ${\mathcal R}^\sharp_{x,a}(W)$ is calculated on the space $\mathscr B_{w} = \mathscr B_{1 +\beta W}$, $\beta := \alpha_0/K$, for some $\alpha_0 \in (0, \alpha)$, and (ii) for all $v \geq u \in \mathscr B_{w}$, 
\begin{align}
 \inf_{x \in B, (x,a)\in \mathbf K} \left\{ \mathcal R_{x,a}(v) - \alpha \nu(v) - \mathcal R_{x,a}(u) +
 \alpha \nu(u) \right\} \geq 0 \label{eq:doeb2}
\end{align}
where $B = \{ x \in \mathbf X: W(x) \leq R \} \in \mathcal
   B(\mathbf X)$ for some $R > 2K/(1-\gamma)$.
\end{assumption}

{\em Comparison with classical MCPs.} In this case, $\mathcal R_{x,a}(v) = \overline{\mathcal R^\sharp}_{x,a}(v) = Q_{x,a}(v)$ and obviously \eqref{eq:doeb2} is equivalent to $\inf_{x \in B, (x,a) \in \mathbf K} Q_{x,a}(A) \geq \nu(A)$ for all $A \in \mathcal B(\mathbf X)$, where $\nu$ becomes a probability measure. Therefore, Assumption \ref{assp:semi} becomes the classical condition with a Lyapunov function, which has been widely used in the MCPs literature (see e.g.~\cite{hernandez1991recurrence,hernandez1999further,vega2003average} and references therein) for studying the optimization problem of average costs with Borel spaces and unbounded costs, to ensure that for each deterministic policy $f \in \Delta_D$, the Markov chain with transition kernel $Q^f$ is $w$-weighted geometric ergodic (for definition see \cite{meyn1993markov}). 

{\em Comparison with entropic maps.} The most widely studied risk map in MCP-literature is the entropic map \eqref{eq:entropic}. While most of them focus on finite or countable state spaces, only a few papers (see e.g.~\cite{di2008infinite} and references therein) discussed the general state space model as we consider in this paper under the sup-norm $\lVert \cdot \rVert_{\infty} := \lVert \cdot \rVert_{w=1}$. Thus we restrict to $\mathscr B$. Furthermore, all of the MCP-literature applying the entropic maps consider merely risk-averse cases, i.e., $\lambda > 0$. Hence, without loss of generality, we set $\lambda = 1$. In the following, we connect the generalized Doeblin's condition stated in \eqref{eq:doeb2} with properties of kernel $Q$. 

\begin{proposition}
 Suppose there exists a set $C \in \mathcal B(\mathbf K)$, a constant $\alpha > 0$ and a probability $\nu$ such that $\inf_{k \in C} Q(\cdot|k) \geq \alpha \nu(\cdot)$. Then for all $v \geq u \in \mathscr B$ and $\lVert u \rVert_\infty \leq M<\infty$, there exists a constant $\epsilon(M) > 0$, such that, for all $v \geq u$, $\lVert u \rVert_\infty \leq M$, 
\begin{equation}
 \inf_{k \in C} \left\{ \mathcal R_{k}(v) - \mathcal R_{k}(u) - \epsilon(M) \nu(v) + \epsilon(M) \nu(u) \right\} \geq 0 \label{eq:en:in}
\end{equation}
where $\mathcal R$ is the entropic risk map defined in \eqref{eq:entropic}. 
\end{proposition}
\begin{proof}
Note that adding any constant valued function to $v$ will not change the inequality in \eqref{eq:en:in}. Hence, we assume that $-2M \leq u \leq 0$. By Lemma 6.1 of \cite{svindland2009subgradients}, we have $g(k,u) := \frac{ \exp(u)}{\int_{\mathbf X} \exp(u) d Q_{k}} \geq \exp(-2M)$ is one subgradient in $\delta \mathcal R_{k}(u)$. Hence, by definition of subgradient and assumption 
\begin{align*}
 \mathcal R_k(v) - \mathcal R_k(u) \geq  \int (v - u) g d Q_{k} \geq \alpha \exp(-2M) \int (v-u) d \nu.
\end{align*}
Setting $\epsilon(M) = \alpha \exp(-2M) > 0$, we obtain \eqref{eq:en:in}. \quad \end{proof}

Note that though \eqref{eq:en:in} is weaker than the generalized Doeblin's condition stated in \eqref{eq:doeb2}, for the purpose of proving the existence of solutions to Poisson equation \eqref{eq:aroe}, it is sufficient (cf.~Proposition 2.2 of \cite{di1999risk}) to consider the functions $u$ upper bounded by some constant $M$ under sup-norm, if the cost function $c$ is upper bounded by $M$, i.e., $\lvert c(k) \rvert \leq M$, $\forall k \in \mathbf K$.  In \cite{di1999risk}, the condition $\sup_{k,k' \in \mathbf K, B \in \mathcal B(\mathbf X)} \lvert Q_k(B) - Q_{k'}(B) \lvert \leq \alpha$ is assumed. It is known that this condition implies (see e.g.~\cite{dobrushin1956central}) the Doeblin's condition holds on the whole space $\mathbf K$, i.e., there exists a probability measure $\mu$ and a constant $\alpha' > 0$ such that $\inf_{k \in \mathbf K} Q(B|k) \geq \alpha' \mu(B), \forall B \in \mathcal B(\mathbf X)$. Hence, setting $W = 1$ in \eqref{eq:ly2}, Assumption \ref{assp:semi} holds. It means that the condition stated in \cite{di1999risk} is covered by our assumption. However, it is still an open question whether our assumption is stronger than the assumption considered in the most recent paper \cite{di2008infinite}.

From now on, we state proofs related to the average risk optimality equation.
\begin{lemma}\label{lm:disc3}
 Suppose Assumption \ref{asmp:sub} and \ref{assp:semi} hold. Then there exists $\bar \alpha \in [0,1)$ such that $\lVert \mathcal F(v) - \mathcal F(u) \rVert_{s,w} \leq \bar \alpha \lVert v - u \rVert_{s,w}$, for all $v,u\in \mathscr B_w$, where $\mathcal F$ is defined in \eqref{eq:F}.
\end{lemma}
\begin{proof}
 Under Assumption \ref{asmp:sub}, by Proposition \ref{prop:sub}, there exists a deterministic policy $f_v, f_u \in \Delta_D$ such that $\mathcal F(v) = \mathcal F^{f_v}(v)$ and $\mathcal F(u) = \mathcal F^{f_u}(u)$. Thus
\begin{align*}
  & \mathcal F(v) - \mathcal F(u) \leq \mathcal F^{f_u}(v) - \mathcal F^{f_u}(u) = \mathcal R^{f_u}(v) - \mathcal R^{f_u}(u) \\
\textrm{and} \quad & \mathcal F(u) - \mathcal F(v) \leq \mathcal F^{f_v}(u) - \mathcal F^{f_v}(v) = \mathcal 
R^{f_v}(u) - \mathcal R^{f_v}(v)
\end{align*}
yield  $\mathcal F_x(v) - \mathcal F_x(u) +  \mathcal F_y(u) - \mathcal F_y(v) \leq \mathcal R^{f_u}_x(v) - \mathcal R^{f_u}_x(u) + \mathcal R^{f_v}_y(u) - \mathcal R^{f_v}_y(v)$ for all  $x,y \in \mathbf X$. By Assumption \ref{assp:semi} and repeating the proof in Theorem \ref{th:semidiscount}, we have
\begin{align*}
 & \mathcal F_x(v) - \mathcal F_x(u) +  \mathcal F_y(u) - \mathcal F_y(v) \\
 \leq & \mathcal R^{f_u}_x(v) - \mathcal R^{f_u}_x(u) + \mathcal R^{f_v}_y(u) - \mathcal R^{f_v}_y(v) \leq \bar \alpha \lVert v - u \rVert_{s,w} d_w(x,y).
\end{align*}
Switching $v$ and $u$, the l.h.s of the inequality does not change. Thus we have $\lvert \mathcal F_x(v) - \mathcal F_x(u) +  \mathcal F_y(u) - \mathcal F_y(v) \rvert \leq \bar \alpha \lVert v - u \rVert_{s,w} d_w(x,y).$
\quad \end{proof}

Above lemma shows that the operator $\mathcal F: \mathscr B_w \rightarrow \mathscr B_w$ is a contraction under the $w$-seminorm. Hence, by  Theorem \ref{th:pe} (i) we immediate obtain the following theorem.
\begin{theorem}\label{th:aroe}
 Suppose Assumption \ref{asmp:sub}, \ref{assmp:disc}(a) and \ref{assp:semi} hold. Then there exists a unique $\rho^* \in \mathbb R$ and $h \in \mathscr B_w$ satisfying the \emph{average risk optimality equation} (AROE)
\begin{align}
 \rho^* + h(x) = \mathcal F(h|x) = \inf_{a\in \mathbf A(x)} \left\{ c(x,a) + \mathcal R(h|x,a) \right\}. \label{eq:aroe}
\end{align}
\end{theorem}

\begin{theorem}\label{th:aroe2}
  Suppose Assumption \ref{asmp:sub}, \ref{assmp:disc}(a) and \ref{assp:semi} hold. Then 
$\rho = J^*(x) = J(x,f^\infty)$ for all $x \in \mathbf X$, where $\rho$ is the solution to the AROE (\eqref{eq:aroe}) and $f$ denotes the optimal selector in the r.h.s of the AROE.
\end{theorem}
\begin{proof}
 Under Assumption \ref{asmp:sub}, \ref{assmp:disc}(a) and \ref{assp:semi}, the existence of one solution to the AROE is guaranteed by Theorem \ref{th:aroe}. By Proposition \ref{prop:sub}, the optimal selector $f$ exists. Thus we assume $\mathcal F(h) = \mathcal F^f (h)$. Due to the AROE,
\begin{align*}
 (\mathcal F^f)^t (h) =  (\mathcal F^f)^{t-1} (\rho + h) = t \rho + h \quad \Rightarrow \quad \lim_{t \rightarrow \infty} \dfrac{1}{t} \lVert (\mathcal F^f)^t (h) - \rho \rVert_w = 0.
\end{align*}
On the other hand, for any $v \in \mathscr B_w$, by Lemma \ref{lm:rhoc} (ii), $\dfrac{1}{t}\lVert (\mathcal F^f)^t (v) - (\mathcal F^f)^t (h) \rVert_w \rightarrow 0$ implies that $J(x,f^\infty) = \lim_{t \rightarrow \infty} \frac{1}{t} (\mathcal F^f)^t_x(0) = \rho, \forall x \in \mathbf X$.

Next we prove that $\rho \leq J(x,\bs \pi)$ for all $\bs \pi \in \Pi_M$ and $x \in \mathbf X$. In fact,  let $\bs \pi = [\pi_0,\pi_1, \ldots]$ be an arbitrary Markov policy. Then by Lemma \ref{lm:rhoc} (ii), for all $v \in \mathscr B_w$,
\begin{align}
 \lim_{t \rightarrow \infty} \frac{1}{t} \lVert (\mathcal F^{\pi_0}(\mathcal F^{\pi_1} \ldots \mathcal F^{\pi_t}(v)) - \mathcal F^{\pi_0}(\mathcal F^{\pi_1} \ldots \mathcal F^{\pi_t}(0))  \rVert_w = 0. \label{eq:tmp5}
\end{align}
By definition $h \leq \mathcal F^\pi(h) - \rho, \forall \pi \in \Delta$. Iterating this inequality yields 
\begin{align}
 h \leq \mathcal F^{\pi_0}( \mathcal F^{\pi_1}( \ldots \mathcal F^{\pi_{t-1}}(h))) - t \rho \Rightarrow \limsup_{t \rightarrow \infty} \frac{1}{t} \mathcal F^{\pi_0}(\mathcal F^{\pi_1} (\ldots \mathcal F^{\pi_{t-1}}(h))) \geq \rho. \label{eq:rhoh}
\end{align}
Note that by definition $J(\bs \pi) = \limsup_{t \rightarrow \infty} \frac{1}{t} \mathcal F^{\pi_0}(\mathcal F^{\pi_1} (\ldots \mathcal F^{\pi_{t-1}}(0)))$ and 
by setting $v=h$ in \eqref{eq:tmp5}, we obtain $J(\bs \pi) = \limsup_{t \rightarrow \infty} \frac{1}{t} \mathcal F^{\pi_0}(\mathcal F^{\pi_1} (\ldots \mathcal F^{\pi_{t-1}}(h)))$. Hence, \eqref{eq:rhoh} implies $\rho \leq J(\bs \pi)$. It follows that $\rho \leq \inf_{\bs \pi \in \Pi_M} J(\bs \pi) = J^*$. Since $f^\infty$ is a valid Markov policy in $\Pi_M$, $\rho = J^* = J(f^\infty)$.
\quad \end{proof}


\section{One Example with Mean-semideviation}
\label{sec:rm}
We present one example by applying the mean-semideviation map defined in \eqref{eq:meansemi} with $r=2$ and $\lambda \in (0,1)$ which has not be considered as an average risk objective in MCP-literature yet. Consider a 1-dimensional simple linear model 
\begin{align}
 X_{t+1} = b(A_t) X_t + N_t \label{eq:dynamic}
\end{align}
where $N_t$ is i.i.d~white noise and $b$ is a real-valued measurable function on $\mathcal B(\mathbf A)$. Hence the transition kernel in this MCP is simply $Q(dy |x,a) = \frac{1}{\sqrt{2\pi}} e^{-\frac{(y-b(a)x)^2}{2}} dy.$ Assume that $\sup_{a \in \mathbf A} \lvert b(a) \rvert = \epsilon < 1$ and set $k = (x,a)$ for simplicity. Recall that $\mathbb E_{k}^Q v := \int_{\mathbf X = \mathbb R} Q(dy|k) v(y)$.

First, we show that $W(x) = x^2$ is a weight function satisfying \eqref{eq:ly2}. Indeed, since the mean-semideviation map $\mathcal R$ is coherent if $\lambda \in (0,1)$, $\overline{\mathcal R^\sharp} = \mathcal R$. 
\begin{align*}
 \mathcal R_{k}(W) = & \mathbb E_{k}^Q W + \lambda \sqrt{ \mathbb E_{k}^Q (W - \mathbb E_{k}^Q W)_+^2}  \leq \mathbb E_{k}^Q W + \lambda \sqrt{ \mathbb E_{k}^Q (W - \mathbb E_{k}^Q W)^2} \\ 
= & b^2(a) x^2 + 1 + \lambda \sqrt{ 2 + 4 b^2(a) x^2} \leq \epsilon^2 x^2 + 1 + \lambda \sqrt{ 2 + 4 \epsilon^2 x^2}
\end{align*}
Note that, for any $\gamma \in (\epsilon^2,1)$, the function $
  t \mapsto (\epsilon^2 - \gamma) t + 1 + \lambda \sqrt{ 2 + 4 \epsilon^2 t}$ 
is concave in $[0,\infty)$ and its maximum is attained at some constant $K$. Thus,
 \begin{align}
 \epsilon^2 x^2 + 1 + \lambda \sqrt{ 2 + 4 \epsilon^2 x^2} \leq \gamma x^2 +K \quad \Rightarrow \quad \mathcal R_{k}(W) \leq \gamma W(x) + K \label{eq:W}
\end{align}
holds for all $k \in \mathbf K$. Thus, the weight function is $w(x) = 1 + \beta W(x) = 1 + \beta x^2$ with some $\beta > 0$.

Next, we check the condition stated in \eqref{eq:doeb2}. Note that if $v \in \mathscr B_w$, then $\mathbb E_{k}^Q [v^2] \leq \lVert v \rVert_w^2  \mathbb E_{k}^Q [1 + \beta y^2]^2 < \infty$ for each $k$. Hence, $L^2(Q(\cdot|k)) \supset \mathscr B_w$. It is known that (see e.g.~\cite{svindland2009subgradients})
\begin{align*}
 g(k,u) = \left\{ \begin{array}{ll}
                   1 & \textrm{if $u$ is constant}\\
                1 - \lambda \dfrac{\mathbb E_{k}^Q \left[ ( u - \mathbb E_{k}^Q u )_+ \right] - ( u - \mathbb E_{k}^Q u )_+}{\sqrt{\mathbb E_{k}^Q [( u - \mathbb E_{k}^Q u )_+^{2}]}} & \textrm{otherwise}
                  \end{array}
 \right.
\end{align*}
is one subgradient defined on $L^2(Q(\cdot|k))$  for $u \in L^2(Q(\cdot|k)) \supset \mathscr B_w$. Note that since
\begin{align*}
 \mathbb E_{k}^Q \left[ ( u - \mathbb E_{k}^Q u )_+ \right] - ( u - \mathbb E_{k}^Q u )_+ \leq \mathbb E_{k}^Q \left[ ( u - \mathbb E_{k}^Q u )_+ \right] \leq \sqrt{\mathbb E_{k}^Q [( u - \mathbb E_{k}^Q u )_+^{2}]},
\end{align*}
$g(k,u) \geq 1 - \lambda > 0$, for $\lambda \in (0,1)$. Besides, it is easy to check that for each $k$ and $u \in \mathscr B_w$,
\begin{align*}
 \mathbb E_{k}^Q [\lvert g(k,u) \rvert w] = \mathbb E_{k}^Q [ g(k,u) w] \leq \mathbb E_{k}^Q \left[ \left(1 + \lambda \frac{( u - \mathbb E_{k}^Q u )_+}{\sqrt{\mathbb E_{k}^Q [( u - \mathbb E_{k}^Q u )_+^{2}]}} \right)w \right]  < \infty.
\end{align*}
Hence $ g(k,u) \in \delta \mathcal R_{k}(u)$, $\forall k \in \mathbf K, u \in \mathscr B_w$. On the other hand, the MCP defined in \eqref{eq:dynamic} satisfies that for any closed set $C_m :\left\{ k \in \mathbf K \mid \lvert x \rvert \leq m \right\}$, $m > 0$, the classical Doeblin's condition holds (cf.~p.~380 of \cite{meyn1993markov}). Hence, by Proposition \ref{prop:doeb}, the generalized Doeblin's condition holds.  Together with \eqref{eq:W}, Assumption \ref{assp:semi} holds. 

\section*{Acknowledgments} The authors would like to thank the reviewers and the associated editor for careful readings and constructive suggestions, which lead to an improvement of the presentation of the paper.

\end{document}